\documentclass[12pt]{article}
by-\headheight \advance \topmargin by-\headsep 
\usepackage{amsmath,amssymb,amsthm,amsfonts,indentfirst}

\newtheorem{thm}{Theorem}[section]

\newtheorem{lem}[thm]{Lemma}
\newtheorem{prop}[thm]{Proposition}

\theoremstyle{definition}

\newtheorem{rem}{Remark}[section]

\numberwithin{equation}{section}

\begin{document}

\newcommand{\la}{\lambda}
\newcommand{\de}{\delta}
\newcommand{\De}{\Delta}
\newcommand{\na}{\nabla}
\newcommand{\Om}{\Omega}
\newcommand{\oo}{\infty}
\newcommand{\os}{\Omega(s)}
\newcommand{\pa}{\partial}
\newcommand{\ti}{\tilde}
\newcommand{\ve}{\varepsilon}
\newcommand{\vf}{\varphi}
\newcommand{\vv}{\bar{V}}
\newcommand{\di}{\na \cdot}
\newcommand{\non}{\nonumber}
\newcommand{\dsf}{\displaystyle \frac}
\newcommand{\Iom}{\displaystyle\int_{\Omega}}
\newcommand{\Ios}{\displaystyle\int_{\Omega(s)}}

\title
{Some Blow-Up Problems For A  Semilinear Parabolic Equation With A
Potential}

\author{  Ting Cheng and Gao-Feng Zheng
\thanks{The corresponding author. }
\\Department of Mathematics,
Huazhong Normal University\\
Wuhan, 430079, P.R. China.\\
tcheng@mail.ccnu.edu.cn, \ \ gfzheng@mail.ccnu.edu.cn}
\date {}
\maketitle

\begin{abstract}
The blow-up rate estimate for the solution to a semilinear parabolic
equation $u_t=\Delta u+V(x) |u|^{p-1}u$  in $\Omega \times (0,T)$
with $0$-Dirichlet boundary condition is obtained. As an
application, it is shown that the asymptotic behavior of blow-up
time and blow-up set of the problem with nonnegative initial data
$u(x,0)=M\vf (x)$ as $M$ goes to infinity, which have been found in
\cite{cer}, are improved under some reasonable and weaker conditions
compared with \cite{cer}.

Key words: Blow-Up rate, Blow-Up time, Blow-Up set, Semilinear
parabolic equations, Potential.

\end{abstract}

\maketitle

\section{Introduction}

In this paper, we are concerned with the following semilinear
parabolic problem
\begin{equation}\label{1.1}
\left\{
\begin{array}{ll}
u_t=\Delta u+V(x) |u|^{p-1}u  & \ \ \mbox{in} \ \ \Omega \times
(0,T),\\
u(x,t)= 0 & \ \ \mbox{on}  \ \partial \Omega \times
(0,T),\\
u(x,0)=u_0(x) & \ \ \mbox{in} \ \ \Omega,
\end{array}
\right.
\end{equation}
where $\Omega \subset \mathbb{R}^N (N\geqslant 3)$ is a bounded,
convex, smooth domain, $1<p<\frac{N+2}{N-2}$, $u_0\in L^{\oo}(\Om)$,
and the potential $V\in C^1(\bar \Om)$ satisfies $V(x)\geqslant c$
for some positive constant $c$ and all $x\in \Om$. It is well-known
that for any $u_0\in L^{\oo}(\Om)$ problem \eqref{1.1} has a unique
local in time solution. Specially, if the $L^{\oo}$-norm of the
initial datum is small enough, then \eqref{1.1} has global,
classical solution, while the solution to \eqref{1.1} ceases to
exist after some time $T>0$ and $\lim\limits_{t\uparrow T}\|u(\cdot
,t)\|_{L^{\oo}(\Om)}=\oo$ provided that the initial datum $u_0$ is
large in some suitable sense. In the latter case we call the
solution $u$ to \eqref{1.1} blows up in finite time and $T$ the
blow-up time. As usual, the blow-up set of the solution $u$ is
defined by
\[
B[u]=\{x\in \bar{\Om}\ | \ \mbox{there exist} \ x_n\to x,
t_n\uparrow T, \mbox{such that} \ |u(x_n,t_n)|\to \oo\}.
\]

Much effort has been devoted to blow-up problems for semilinear
parabolic equations since the pioneering works in 1960s due in
particular to interest in understanding the mechanism of thermal
runaway in combustion theory and as a model for reaction-diffusion.
See for example, \cite{ba,be,ch,fm,f,gv,k,l}. The seminal works to
problem \eqref{1.1} with $V(x)\equiv 1$ were done by Giga-Kohn
\cite{gk1,gk2,gk3}. In their paper \cite{gk2}, among other things,
they have obtained a blow-up rate estimate, which is crucial to
obtain the asymptotic behavior of the blow-up solution near the
blow-up time. More precisely, under the assumptions that the domain
$\Om$ is the entire space or convex and the solution is nonnegative
or $1<p<\frac{3N+8}{3N-4} (N\geqslant 2)$ or $1<p<\oo (N=1)$, they
proved that
\[
|u(x,t)|\leqslant C(T-t)^{\frac{1}{p-1}}, \ \ \forall \ (x,t)\in \Om
\times (0,T),
\]
where $C>0$ is a constant and $T>0$ is the blow-up time. More
recently, the same estimate has been obtained by
Giga-Matsui-Sasayama \cite{gms1, gms2} for any subcritical $p$
(i.e., $1<p<\frac{N+2}{N-2}$ when $N\geqslant 3$, $1<p<\oo$ when
$N=1,2$).

Whether the similar blow-up rate estimate holds for the problem
\eqref{1.1} for general potential $V$, to our best knowledge, is not
well-understood up to now. Our first goal in this paper is to give
an affirmative answer to this question. We have the following

\begin{thm}\label{thm1}
Let $u$ be a blow-up solution to \eqref{1.1} with a blow-up time
$T$. There exists a positive constant $C$ depending only on
$n,p,\Om$, a bound for $T^{1/(p-1)}\|u_0\|_{L^{\oo}(\Om)}$ and the
positive lower bound $c$ for $V$ and $\|V\|_{C^1(\bar\Om)}$, such
that
\begin{equation}\label{1.2}
\|u(\cdot, t)\|_{L^{\oo}(\Om)}\leqslant C(T-t)^{-1/(p-1)}, \ \ \
\forall \ t\in (0,T).
\end{equation}
\end{thm}

As in \cite{gk2}, we convert our problem to a uniform bound for a
global in time solution $w$ of the rescaled equation
\[
 w_s-\De w+\frac{1}{2}y\cdot \na w +\beta w -\bar V |w|^{p-1}w
 =0,  \ \ \ \beta=\frac{1}{p-1},
\]
with
\[
w(y,s)=(T-t)^\beta u(a+y\sqrt{T-t},t), \ \bar V(y,s)=V(a+ye^{-s/2}),
\]
where $a\in \Om$ is the center of the rescaling.

The proof of Theorem \ref{thm1} depends heavily on the methods
developed by Giga-Kohn in \cite{gk2} and Giga-Matsui-Sasayama in
\cite{gms1,gms2}. However our result is definitely not a direct
consequence of their works. Due to the appearance of the potential
$V$, some extra works should be done. It turns out that the key
point and the main difference is to establish an upper bound for the
global energy of $w$ given by
\[
E[w](s)=\dsf{1}{2}\int_{\os} (|\na w|^2 +\beta w^2)\rho \,dy
     -\dsf{1}{p+1}\int_{\os} \vv |w|^{p+1}\rho \,dy,
\]
where $\rho(y)=e^{-\frac{|y|^2}{4}}$. A lower bound for the energy
can be obtained without much effort. When $V \equiv 1$, these bounds
come easily from the Liapunov structure of the equation, i.e., the
energy $E[w]$ is non-increasing in time. In our case this does not
hold anymore. There is a ``bad" term
\[
\int_{\os}\frac{\partial \bar V}{\partial s}|w|^{p+1}\rho \,dy
\]
involved in the derivative of the energy $E[w]$. The main idea of
proving Theorem \ref{thm1} is as follows: First, using the fact that
$\frac{\partial \bar V}{\partial s}$ is uniformly bounded in $\os$
for all $s$, we get a rough control of the growth of the global
energy $E[w]$. Since the term $\frac{\partial \bar V}{\partial s}$
can be written as $\na V(x)\cdot y e^{-s/2}$, we can use the
information of the decay term $e^{-s/2}$. However, it has
disadvantage that the unbounded thing $y$ involves. Therefore we
need some information from \textit{higher level energies}
\[
E_{2k}[w](s)=\dsf{1}{2}\int_{\os} (|\na w|^2 +\beta w^2)|y|^{2k}\rho
\,dy -\dsf{1}{p+1}\int_{\os} \vv |w|^{p+1}|y|^{2k}\rho \,dy, \ \ \
k\in \mathbb{N}.
\]
So our second step is to establish the control of the higher level
energies. An upper bound for $E_{2k}[w]$ is obtained by an integral
involving lower level energy. A lower bound for $E_{2k}[w]$ is
obtained by two inequalities involving $\dsf{d}{ds}\int_{\os}
w^2|y|^{2k}\rho \,dy$ and $dE_{2k}[w]/ds$. Finally we obtain an
upper bound for the energy $E[w]$. To this end, the growth of lower
level energies is improved by applying the growth of the higher
level energies. An upper bound of the global energy $E[w]$ is
obtained by a similar trick to bootstrap argument. Once these bounds
are in hands, similar arguments to \cite{gms1,gms2} can be applied
to show the boundedness of the global in time solution $w$, which in
turn implies the blow-up rate estimate \eqref{1.2}.

Another aim of this paper is to establish the asymptotic behavior of
blow-up time and blow-up set of the blow-up solution to the problem
\eqref{1.1} with nonnegative initial data $u_0=M \vf$ as $M\to \oo$.
In this case, the problem we focused on can be rewritten as
\begin{equation}\label{1.3}
\left\{
\begin{array}{ll}
u_t=\Delta u+V(x) u^{p}  & \ \ \mbox{in} \ \ \Omega \times
(0,T),\\
u(x,t)= 0 & \ \ \mbox{on}  \ \partial \Omega \times
(0,T),\\
u(x,0)=M\vf(x) & \ \ \mbox{in} \ \ \Omega,
\end{array}
\right.
\end{equation}
where $\vf \in C(\bar \Om)$ satisfies $\vf|_{\partial \Om} =0,
\vf(x)>0, \ \forall \ x\in \Om$ and $V$ satisfies the same
conditions as before. For these issues of blow-up problems to
\eqref{1.3}, we improve the results which have been obtained by
Cortazar-Elgueta-Rossi \cite{cer} recently.

In \cite{cer}, they have made some more technical condition on
$\vf$:
\begin{equation}\label{1.5}
M\Delta \vf +\frac{1}{2}\min_{x\in \Om}V(x)M^p\vf^p \geqslant 0.
\end{equation}
The assumptions on $\Om, p$ and $V$ are the same as ours (although
their assumption that $V$ is Lipschitz is replaced by $V\in C^1(\bar
\Om)$ in our case, our results still hold when $V$ is Lipschitz).
Under these assumptions, they proved that there exists $\bar M>0$
such that if $M>\bar M$, then blow-up occurs and the blow-up time
$T(M)$ and the blow-up set $B[u]$ of the blow-up solution to
\eqref{1.3} satisfy
\[
-\frac{C_1}{M^{\frac{p-1}{4}}}\leqslant
T(M)M^{p-1}-\frac{A}{p-1}\leqslant \frac{C_2}{M^{\frac{p-1}{3}}},
\]
\[
\vf^{p-1}(a)V(a)\geqslant \frac{1}{A}-\frac{C}{M^{\gamma}}, \ \
\mbox{ for all } \ a\in B[u],
\]
where $A=\left(\max_{x\in \Om}\vf^{p-1}(x)V(x)\right)^{-1}, \
\gamma=\min(\frac{p-1}{4}, \frac{1}{3})$ and $C_1, C_2$ are two
positive constants.

For the upper bound estimate on blow-up time, we have the following
\begin{thm}\label{thm2}
Let $\Om\subset \mathbb{R}^N (N\geqslant 3)$ be a smooth bounded
domain, $p>1$, $V, \vf$ be continuous functions on $\bar \Om$ with
$\vf|_{\partial \Om}=0, \vf(x)>0, V(x)\geqslant c, \ \forall \ x\in
\Om$ for some $c>0$. Then for any $k>p-1$ there exists a constant
$C>0$ and $M_0>0$ such that for every $M\geqslant M_0$, the solution
to \eqref{1.3} blows up in finite time that verifies
\begin{equation}\label{1.4}
T(M)\leqslant \frac{A}{(p-1)M^{p-1}}+ C M^{-k},
\end{equation}
where $A=\left(\max_{x\in \Om}\vf^{p-1}(x)V(x)\right)^{-1}$.
\end{thm}
\begin{rem}
Our assumptions are weaker than ones in \cite{cer}. In \cite{cer},
they required $V$ and $\vf$ are Lipschitz continuous. Furthermore,
our result tells that the decay of the upper bound of
$T(M)-\frac{A}{(p-1)M^{p-1}}$ can be faster than which has been
obtained in \cite{cer}.
\end{rem}

Notice that the proof of the upper bound of blow-up time in
\cite{cer} depends on an argument of so-called ``projection method"
(see e.g. \cite{k}) and the essential assumption that $V, \vf$ are
Lipschitz continuous. Our proof of Theorem \ref{thm2} requires a
$L^2$-method (see e.g. \cite{ba}). The advantage of this method
compared with one in \cite{cer} is that we do not need to control
the first eigenvalue of Laplacian with Dirichlet boundary condition.

For the lower bound estimate for the blow-up time and the asymptotic
behavior of blow-up set, we have
\begin{thm}\label{thm3}
Let $\Om\subset \mathbb{R}^N (N\geqslant 3)$ be a convex, bounded,
smooth domain, $1<p<\frac{N+2}{N-2}$, $\vf$ be a continuous function
on $\bar \Om$ with $\vf|_{\partial \Om}=0, \vf(x)>0, \ \forall \
x\in \Om$, and $V\in C^1(\bar \Om)$ with $V(x)>c, \ \forall \ x\in
\Om$ for some $c>0$. Then there exist two positive constants $C_1,
C_2$ such that
\begin{equation}\label{1.6}
T(M)M^{p-1}\geqslant -\frac{C_1}{M^{\frac{p-1}{4}}}
\end{equation}
\begin{equation}\label{1.7}
\vf^{p-1}(a)V(a)\geqslant \frac{1}{A}-\frac{C_2}{M^{\frac{p-1}{4}}},
\ \ \mbox{ for all } \ a\in B[u],
\end{equation}
where $A=\left(\max_{x\in \Om}\vf^{p-1}(x)V(x)\right)^{-1}$.
\end{thm}

Applying Theorem \ref{thm1} and the method in \cite{cer}, we get
Theorem \ref{thm3} immediately. The only difference is that the role
of Lemma 2.1 in \cite{cer} is replaced by that of our Theorem
\ref{thm1} now.

\begin{rem}
In our case, we do not need the assumption \eqref{1.5} anymore.
\end{rem}

\begin{rem}
As described in \cite{cer}, the asymptotics depend on a combination
of the shape of both $\vf$ and $V$. To see this, if we drop the
Laplacian, we get the ODE $u_t=V(x)u^p$ with initial condition
$u(x,0)=M\vf(x)$. This gives $u(x,t)=C(T-t)^{-1/(p-1)}$ with
\[
T=\frac{M^{1-p}}{(p-1)V(x)\vf^{p-1}(x)}.
\]
It turns out that blow-up occurs at point $x_0$ such that
$V(x_0)\vf^{p-1}(x_0)=\max\limits_{x\in \Om}V(x)\vf^{p-1}(x)$. So
the quantity $\max\limits_{x\in \Om}V(x)\vf^{p-1}(x)$ plays a
crucial role in the problem.
\end{rem}

\begin{rem}
Also as in \cite{cer}, \eqref{1.7} shows that the blow-up set
concentrates when $M\to \oo$ near the set where $\vf^{p-1}V$ attains
its maximum. Notice that $1/A=\vf^{p-1}(\bar a)V(\bar a)$ for any
maximizer $\bar a$. If $\bar a$ is a non-degenerate maximizer, we
conclude that there exist constants $c,d>0$ such that
\[
\vf^{p-1}(\bar a)V(\bar a)-\vf^{p-1}(x)V(x)\geqslant c|\bar a-x|^2 \
\ \ \mbox{for all } \ x\in B(\bar a, d).
\]
So \eqref{1.7} implies
\[
|\bar a-a|\leqslant \frac{C}{M^{(p-1)/8}}, \ \ \forall \ a\in B[u].
\]
\end{rem}

Throughout the paper we will denote by $C$ a constant that does not
depends on the solution itself. And it may change from line to line.
And $K_1,K_2,\cdots$, $L_1,L_2,\cdots$, $M_1, M_2,\cdots$,
$N_1,N_2,\cdots$, $Q_1, Q_2,\cdots$ are positive constants depending
on $p, N, \Om$, a lower bound of $V$, $\|V\|_{C^1(\bar \Om)}$ and
the initial energy
$E[w_0]$. Here and hereafter $w_0(y)=w(y,s_0)$.\\

\noindent{\bf Acknowledgment.} This work is partially supported by
NSFC No.10571069.

\section{Blow-Up Rate Estimates}

In this section, we will prove Theorem \ref{thm1}.

We introduce the rescaled function
\begin{equation}\label{3.1}
 w^{a}(y,s)=(T-t)^\beta u(a+y\sqrt{T-t},t)
\end{equation}
with $s=-\log(T-t), \ \beta=\frac{1}{p-1}.$ We shall denote $w^a$ by
$w$. If $u$ solves \eqref{1.1}, then $w$ satisfies
\begin{equation}\label{3.2}
 w_s-\De w+\frac{1}{2}y\cdot \na w +\beta w -|w|^{p-1}w V(a+ye^{-s/2})=0
  \ \ \ \mbox{\rm in} \ \os \times (s_0, \oo)
\end{equation}
where $\os=\Om_a(s)=\{y: a+y e^{-s/2} \in \Om\}, \  s_0=-\log T.$

We may assume $T=1$ as in \cite{gms1} so that we assume $s_0=0.$
Here and hereafter we may denote $V(a+ye^{-s/2})$ by $\vv(y,s).$

By introducing a weight function
$\rho(y)=\exp\left(-\dsf{|y|^2}{4}\right)$,  we can rewrite
\eqref{3.2} as the divergence form:
\begin{equation}\label{3.4}
 \rho w_s=\na \cdot (\rho \na w) -\beta \rho w +\vv |w|^{p-1}w \rho
  \ \ \ \mbox{\rm in} \ \os \times (0, \oo).
\end{equation}
As stated in \cite{gms1}, we may assume
$$w, w_s, \na w \ \mbox{and} \ \na^2 w \ \mbox{are bounded and continuous on}
  \ \os \times [0, s] \  \mbox{for all} \ s<\oo.$$

\subsection{Global energy estimates}

We introduce the energy of $w$ of the form (we call it the ``global
energy")
$$E[w](s)=\dsf{1}{2}\int_{\os} (|\na w|^2 +\beta w^2)\rho \,dy
     -\dsf{1}{p+1}\int_{\os} \vv |w|^{p+1}\rho \,dy.$$
We shall show that this global energy satisfies the following
estimates.
\begin{prop}\label{p3.1}
 Let $w$ be a global solution of \eqref{3.4}, then
 \begin{equation}\label{3.5}
  -K_1 \leqslant E[w] \leqslant K_2.
 \end{equation}
\end{prop}
\begin{prop}\label{p3.2}
 Let $w$ be a global solution of \eqref{3.4}, then
 \begin{eqnarray}
  \int_0^\oo \|w_s;L_\rho^2(\os)\|^2 ds \leqslant N_1, \label{3.6}\\
  \|w;L_\rho^2(\os)\|^2 \leqslant N_2, \label{3.7}\\
  \int_s^{s+1}\|w;L_\rho^{p+1}(\os)\|^{2(p+1)} ds \leqslant N_3. \label{3.8}
 \end{eqnarray}
\end{prop}
We will prove these two properties in the following subsections.

\subsubsection{Lower bound for $E[w]$}

\begin{lem}\label{lem}
 $E[w] \geqslant -K_1.$
\end{lem}

We see from \eqref{3.4} that
\begin{eqnarray}\label{3.9}
 \frac{1}{2}\frac{d}{ds}\int_{\os} w^2 \rho \,dy
 &=& \int_{\os} w w_s \rho \,dy
    =\int_{\os} w(\na \cdot (\rho \na w) -\beta \rho w +\vv |w|^{p-1}w \rho) \,dy\non\\
 &=& -\int_{\os} |\na w|^2 \rho \,dy
     -\int_{\os} \beta w^2 \rho \,dy
     +\int_{\os} \vv |w|^{p+1} \rho \,dy \non\\
 &=& -2E[w] +\frac{p-1}{p+1}\int_{\os} \vv |w|^{p+1} \rho \,dy.
\end{eqnarray}
Calculating the derivative of $E[w]$ and noting that $w_s|_{\pa
\os}=-\frac{1}{2}y \cdot \na w$ we have
\begin{eqnarray}\label{3.10}
 \frac{d}{ds}E[w](s)
 &=& \int_{\os} (\na w \cdot \na w_s +\beta w w_s)\rho \,dy
     -\int_{\os} \vv |w|^{p-1}w w_s \rho \,dy \non\\
 & & +\frac{1}{4}\int_{\pa\os} |\na w|^2(y \cdot \gamma)\rho \,d\sigma
     -\frac{1}{p+1}\int_{\os} \frac{\pa \vv}{\pa s}|w|^{p+1}\rho \,dy \non\\
 &=& -\int_{\os} \na \cdot (\rho \na w)w_s \,dy
     +\int_{\pa\os} (\rho \na w \cdot \gamma)w_s \,d\sigma
     +\int_{\os} \beta w w_s \rho  \,dy \non\\
 & & -\int_{\os} \vv |w|^{p-1}w w_s \rho  \,dy
     +\frac{1}{4}\int_{\pa\os} |\na w|^2(y \cdot \gamma)\rho \,d\sigma
     -\frac{1}{p+1}\int_{\os} \frac{\pa \vv}{\pa s}|w|^{p+1}\rho \,dy \non\\
 &=& -\int_{\os} \na \cdot (\rho \na w)w_s  \,dy
     +\int_{\os} \beta w w_s \rho \,dy
     -\int_{\os} \vv |w|^{p-1}w w_s \rho \,dy \non\\
 & & -\frac{1}{4}\int_{\pa\os} |\na w|^2(y \cdot \gamma)\rho \,d\sigma
     +\frac{1}{2(p+1)}\int_{\os} \na \vv \cdot y |w|^{p+1}\rho \,dy \non\\
 &=& -\int_{\os} w_s^2 \rho \,dy
     -\frac{1}{4}\int_{\pa\os} |\na w|^2(y \cdot \gamma)\rho \,d\sigma \non\\
 & & +\frac{1}{2(p+1)}\int_{\os} \na \vv \cdot y |w|^{p+1}\rho \,dy
\end{eqnarray}
or
\begin{eqnarray}\label{3.11}
 \int_{\os} w_s^2 \rho \,dy
 &=& -\frac{d}{ds}E[w](s)
     -\frac{1}{4}\int_{\pa\os} |\na w|^2(y \cdot \gamma)\rho \,d\sigma \non\\
 & & +\frac{1}{2(p+1)}\int_{\os} \na \vv \cdot y |w|^{p+1}\rho \,dy.
\end{eqnarray}
Notice that $\vv$ is bounded. By \eqref{3.9}, using Young's
inequality, we have
\begin{eqnarray*}
 -2E[w]+C\int_{\os} |w|^{p+1}\rho \,dy
 & \leqslant & -2E[w]+\frac{p-1}{p+1}\int_{\os} \vv |w|^{p+1}\rho \,dy \\
 &  =  & \int_{\os} w w_s \rho \,dy \\
 & \leqslant & \ve \int_{\os} w_s^2 \rho \,dy
         +\ve \int_{\os} |w|^{p+1} \rho \,dy
         +C(\ve).
\end{eqnarray*}
Taking $\ve$ small enough we get
\begin{equation}\label{3.12}
 \int_{\os} |w|^{p+1}\rho \,dy
  \leqslant  C E[w] +\ve \int_{\os} w_s^2 \rho \,dy
          +C(\ve).
\end{equation}
Since $\sup\limits_{y \in \os}|\na \vv||y|
       =\sup\limits_{x \in \Om}|\na V||x-a|$
is bounded and $\Om$ is convex, it follows from \eqref{3.10} and
\eqref{3.12} that
\begin{eqnarray*}
 \frac{d}{ds}E[w](s)
 & \leqslant & -\int_{\os} w_s^2 \rho \,dy
         +C \int_{\os} |w|^{p+1} \rho \,dy  \\
 & \leqslant & -(1-\ve)\int_{\os} w_s^2 \rho \,dy
         +C E[w] +C(\ve).
\end{eqnarray*}
Take $\ve$ small then we have
\begin{equation}\label{3.13}
 \frac{d}{ds}E[w](s) \leqslant C_1 E[w] +C_2.
\end{equation}
From this inequality, we claim that $E[w] \geqslant
-\frac{C_2}{C_1}.$ If not, then there exists $s_1>0$ such that
$E[w](s_1)<-\frac{C_2}{C_1}$. By \eqref{3.13}, we have
$\frac{d}{ds}E[w](s_1)<0.$ This implies that
$$E[w](s)<-\frac{C_2}{C_1} \ \ \mbox{\rm for all} \ s \geqslant s_1.$$
Hence by \eqref{3.9} and Jensen's inequality, for $s \geqslant s_1,$
we have
$$\frac{1}{2}\frac{d}{ds}\int_{\os} w^2 \rho \,dy
  \geqslant C \int_{\os} |w|^{p+1} \rho \,dy
  \geqslant C\left(\int_{\os} w^2 \rho \,dy\right)^{\frac{p+1}{2}}.$$
This fact shows that $\int_{\os} w^2 \rho \,dy$ will blow up in
finite time, which is impossible.

\subsubsection{Upper bound for $E[w]$}

To find an upper bound for $E[w]$, we introduce
$$E_{2k}[w]=\dsf{1}{2}\int_{\os} (|\na w|^2 +\beta w^2)|y|^{2k} \rho \,dy
        -\dsf{1}{p+1}\int_{\os} \vv |w|^{p+1} |y|^{2k} \rho \,dy,
        \ \ k\in \mathbb{N}.$$
For this energy functional, we shall prove the following properties.
\begin{prop}\label{p3.3}
 \begin{eqnarray}\label{3.14}
  \frac{1}{2}\frac{d}{ds}\int_{\os} w^2 \rho |y|^{2k}\,dy
  &=& -2E_{2k}[w]+\frac{p-1}{p+1}\int_{\os} \vv |w|^{p+1} \rho |y|^{2k}\,dy \non\\
  & & +\int_{\os} k\left(n+2k-2-\frac{1}{2}|y|^2\right)w^2 |y|^{2k-2}\rho \,dy.
 \end{eqnarray}
\end{prop}
\begin{prop}\label{p3.4}
 \begin{eqnarray}\label{3.15}
  \int_{\os} w_s^2 \rho |y|^{2k}\,dy
  &=& -\frac{d}{ds}E_{2k}[w]
      -2k\int_{\os} \rho (y \cdot \na w)w_s |y|^{2k-2}\,dy
      -\frac{1}{4}\int_{\pa\os} \left|\frac{\pa w}{\pa \gamma}\right|^2
      (y \cdot \gamma)\rho |y|^{2k}\,d\sigma \non\\
  & & -\frac{1}{p+1}\int_{\os} \frac{\pa \vv}{\pa s} |w|^{p+1} \rho |y|^{2k}\,dy.
 \end{eqnarray}
\end{prop}
\textbf{Proof of Proposition \ref{p3.3}.} Similar to that of
\cite{gk2} Proposition 4.1.

\textbf{Proof of Proposition \ref{p3.4}.}
\begin{eqnarray*}
 \frac{d}{ds}E_{2k}[w]
 &=& \int_{\os} (\na w \cdot \na w_s +\beta w w_s -\vv |w|^{p-1}w w_s)\rho |y|^{2k} \,dy \\
 & & -\frac{1}{p+1}\int_{\os} \frac{\pa \vv}{\pa s}|w|^{p+1}\rho|y|^{2k} \,dy
     +\frac{1}{4}\int_{\pa\os} |\na w|^2(y \cdot \gamma)\rho |y|^{2k} \,d\sigma.
\end{eqnarray*}
Estimating the first term of the right hand side, we get
\begin{eqnarray*}
 \int_{\os} \na w \cdot \na w_s \rho |y|^{2k} \,dy
 &=& -\int_{\os} \na \cdot (\rho |y|^{2k} \na w)w_s \,dy
     +\int_{\pa\os} \rho |y|^{2k} \na w \cdot \gamma w_s \,d\sigma \\
 &=& -\int_{\os} \na \cdot (\rho \na w)w_s |y|^{2k} \,dy
     -2k\int_{\os} w_s \rho \na w \cdot y  |y|^{2k-2}\,dy\\
 & & -\frac{1}{2}\int_{\pa\os} |\na w|^2(y \cdot \gamma)\rho |y|^{2k} \,d\sigma.
\end{eqnarray*}
Hence we have
\begin{eqnarray*}
 \frac{d}{ds}E_{2k}[w]
 &=& -\int_{\os} w_s\left(\na \cdot (\rho \na w) +\beta w \rho -\vv w^p \rho\right) |y|^{2k} \,dy
     -2k\int_{\os} w_s \rho \na w \cdot y |y|^{2k-2} \,dy\\
 & & -\frac{1}{p+1}\int_{\os} \frac{\pa \vv}{\pa s}|w|^{p+1}\rho|y|^{2k} \,dy
     -\frac{1}{4}\int_{\pa\os} |\na w|^2(y \cdot \gamma)\rho |y|^{2k} \,d\sigma  \,dy\\
 &=& -\int_{\os} w_s^2 \rho |y|^{2k}\,dy
     -2k\int_{\os} w_s \rho \na w \cdot y |y|^{2k-2} \,dy \\
 & &   -\frac{1}{p+1}\int_{\os} \frac{\pa \vv}{\pa s}|w|^{p+1}\rho|y|^{2k} \,dy
     -\frac{1}{4}\int_{\pa\os} |\na w|^2(y \cdot \gamma)\rho |y|^{2k} \,d\sigma.
\end{eqnarray*}

For $k=1$, similar to Proposition 4.2 of \cite{gk2} we now state an
parabolic type Pohozaev identity.

\begin{prop}\label{p3.5}
 \begin{eqnarray}\label{3.16}
  & & \frac{1}{2}\frac{d}{ds}\int_{\os} \left(\frac{1}{2}|y|^2-n\right)w^2 \rho \,dy
      -(p+1)\int_{\os} (y \cdot \na w)w_s \rho \,dy \non\\
  &=& \int_{\os} |\na w|^2 \rho \left(c_2+\frac{p-1}{4}|y|^2\right)\,dy
      -\frac{p+1}{2}\int_{\pa\os} \left|\frac{\pa w}{\pa \gamma}\right|^2(y \cdot \gamma)\rho \,d\sigma \non\\
  & & +\int_{\os} \na \vv \cdot y |w|^{p+1} \rho \,dy.
 \end{eqnarray}
\end{prop}

We now define
\begin{equation}\label{3.19}
 \ti E_2[w]
 \triangleq
 E_2[w]-\frac{1}{2}\int_{\os} \left(\frac{1}{2}|y|^2-n\right)w^2 \rho \,dy.
\end{equation}

\begin{lem}\label{lem3.7}
 \begin{equation}\label{3.24}
  \frac{d(\ti E_2 +c_3 E)}{ds}
  \leq
  -c_4\int_{\os}(w_s^2+|\na w|^2)(1+|y|^2) \rho \,dy
  +\lambda (\ti E_2 +c_3 E) +c_5,
 \end{equation}
 where $\lambda=\dsf{8}{p-1}\dsf{d_2}{d_1}$ and $c_5$ depends on $p, d_1, d_2, \eta,$ $d_1$ and $d_2$
 are constants such that $ V(x) \geqslant d_1 >0$ and
 $\sup\limits_{x\in \Om}|\na V(x)| \mbox{diam} (\Om) \leqslant 2d_2$ and $\eta$ is a small constant.
\end{lem}
 \textbf{Proof.} By \eqref{3.15} and \eqref{3.16} we obtain that
  \begin{eqnarray}\label{3.20}
   \frac{d\ti E_2}{ds}
   &=& -\int_{\os}|w_s|^2 \rho |y|^2 \,dy
       -(p+3)\int_{\os} (y \cdot \na w)w_s \rho \,dy
       -\frac{1}{4}\int_{\pa \os} (y \cdot \gamma)
        \left|\frac{\pa w}{\pa \gamma}\right|^2 \rho |y|^2 \,d\sigma \non\\
   & & -\int_{\os} |\na w|^2\left(c_2+\frac{p-1}{4}|y|^2\right) \rho \,dy
       +\frac{p+1}{2} \int_{\pa \os} (y \cdot \gamma)
        \left|\frac{\pa w}{\pa \gamma}\right|^2 \rho \,d\sigma \non\\
   & & -\frac{1}{p+1} \int_{\os} \frac{\pa \vv}{\pa s} |w|^{p+1} \rho |y|^2 \,dy
       +2\int_{\os} \frac{\pa \vv}{\pa s} |w|^{p+1} \rho \,dy.
  \end{eqnarray}
  Since $\Om$ is convex, the third term on the right is always negative.
  We control the second term by applying the Cauchy-Schwarz inequality:
  for any $\ve >0$,
  $$\left|\int_{\os} (y \cdot \na w)w_s \rho \,dy\right|
    \leqslant \ve \int_{\os} \rho |y|^2 |\na w|^2 \,dy
        +\frac{1}{4 \ve} \int_{\os} \rho |w_s|^2 \,dy.$$
  Choosing $\ve$ small enough that
  $\frac{p-1}{4}-(p+3)\ve =\de >0,$
  we conclude that
  \begin{eqnarray*}
   \frac{d\ti E_2}{ds}
   &\leqslant& -\int_{\os} (|w_s|^2|y|^2+\de |\na w|^2|y|^2+c_2|\na w|^2)\rho\,dy\\
   &   & +\frac{p+1}{2}\int_{\pa \os} (y \cdot \gamma)
          \left|\frac{\pa w}{\pa \gamma}\right|^2 \rho \,d\sigma
         +\frac{p+3}{4\ve}\int_{\os} \rho |w_s|^2 \,dy\\
   &   & -\frac{1}{p+1} \int_{\os} \frac{\pa \vv}{\pa s} |w|^{p+1} \rho |y|^2 \,dy
         +2 \int_{\os} \frac{\pa \vv}{\pa s} |w|^{p+1} \rho \,dy.
  \end{eqnarray*}
  Now choose $c_3 > \max (2(p+1), 1+\frac{p+3}{4 \ve})$, and apply \eqref{3.11} to get
  $$\frac{p+1}{2}\int_{\pa \os} (y \cdot \gamma)
     \left|\frac{\pa w}{\pa \gamma}\right|^2 \rho \,d\sigma
    +\left(1+\frac{p+3}{4 \ve}\right)\int_{\os} \rho |w_s|^2 \,dy
    +c_3\frac{dE}{ds}
    \leqslant
    -\frac{c_3}{p+1}\int_{\os} \frac{\pa \vv}{\pa s} |w|^{p+1} \rho \,dy.$$
  Let $2c_4=\min(1, \de, c_2) > 0$, we derive that
  \begin{eqnarray}\label{3.21}
   \frac{d(\ti E_2 +c_3 E)}{ds}
   &\leq&
    -2c_4\int_{\os}(w_s^2+|\na w|^2)(1+|y|^2) \rho \,dy
    +2 \int_{\os} \frac{\pa \vv}{\pa s} |w|^{p+1} \rho \,dy\non\\
   &    &
    -\frac{c_3}{p+1}\int_{\os} \frac{\pa \vv}{\pa s} |w|^{p+1} \rho \,dy
    -\frac{1}{p+1}\int_{\os} \frac{\pa \vv}{\pa s} |w|^{p+1} \rho |y|^2 \,dy \\
   & \leqslant&
   -2c_4\int_{\os}(w_s^2+|\na w|^2)(1+|y|^2) \rho \,dy
   +\frac{2}{p+1}\int_{\os} (c_3 +|y|^2)\frac{\pa \vv}{\pa s} |w|^{p+1} \rho \,dy.  \non
  \end{eqnarray}
  Note that $ \vv(y,s) \geqslant d_1 >0.$
  From \eqref{3.9} we get
  $$\frac{p-1}{p+1}d_1\int_{\os} |w|^{p+1} \rho \,dy
    \leqslant
    \frac{p-1}{p+1}\int_{\os} \vv |w|^{p+1} \rho \,dy
    =
    2E[w]+\int_{\os} w w_s \rho \,dy.$$
  In the following we will denote $\frac{p+1}{(p-1)d_1}$ by $c(p,d_1)$.
  Making use of the inequality
  \begin{equation}\label{*}
   ab \leqslant \ve(a^2+b^{p+1})+C(\ve), \ \ p>1, \ \forall \ \ve >0,
  \end{equation}
  we obtain that
  \begin{eqnarray*}
   \int_{\os} |w|^{p+1} \rho \,dy
   &\leqslant& 2c(p,d_1)E[w] +\int_{\os} w w_s c(p,d_1) \rho \,dy \\
   &\leqslant& 2c(p,d_1)E[w] +\eta\int_{\os} w^{p+1} \rho \,dy
         +\eta \int_{\os} w_s^{2} \rho \,dy
         +C(p,d_1,\eta).
  \end{eqnarray*}
  Here and hereafter $C(p,d_1,\eta)$ denotes a constant depending on
  $p,d_1,\eta$ and may be different at each occurrence. Take $\eta <1$
  and we hence have
  \begin{equation}\label{3.22}
   \int_{\os} w^{p+1} \rho \,dy
   \leqslant  \frac{2c(p,d_1)}{1-\eta}E[w]
         +\frac{\eta}{1-\eta} \int_{\os} w_s^{2} \rho \,dy
         +C(p,d_1,\eta).
  \end{equation}
  From \eqref{3.14} we obtain that
  \begin{eqnarray*}
   \frac{p-1}{p+1}d_1\int_{\os} |w|^{p+1} \rho |y|^2 \,dy
    &\leqslant&
    \frac{p-1}{p+1}\int_{\os} \vv |w|^{p+1} \rho |y|^2 \,dy \\
    & = &
      2E_2[w]+\int_{\os} w w_s \rho |y|^2 \,dy
       -\int_{\os} \left(n-\frac{1}{2}|y|^2\right)w^2 \rho \,dy \\
    &\leqslant&
     2\ti E_2[w]+\int_{\os} |w w_s| |y|^2  \rho \,dy
       +2\int_{\os} \left(\frac{1}{2}|y|^2-n\right)w^2 \rho \,dy.
  \end{eqnarray*}
  Thanks to \eqref{*}, we hence get
  \begin{eqnarray*}
   \int_{\os} |w|^{p+1} \rho |y|^2 \,dy
   &\leqslant& 2c(p,d_1)\ti E_2[w]
         +\int_{\os} w^2|y|^{\frac{4}{p+1}}
           \cdot c(p,d_1)|y|^{\frac{2(p-1)}{p+1}} \cdot \rho \,dy \\
   &   & +\frac{\eta}{2} \int_{\os} |w|^{p+1} \rho |y|^2 \,dy
         +\frac{\eta}{2} \int_{\os} w_s^{2} \rho |y|^2 \,dy
         +C(p,d_1,\eta) \\
   &\leqslant& 2c(p,d_1)\ti E_2[w]
         +\eta \int_{\os} |w|^{p+1} \rho |y|^2 \,dy
         +\frac{\eta}{2} \int_{\os} w_s^{2} \rho |y|^2 \,dy
         +C(p,d_1,\eta).
  \end{eqnarray*}
  Therefore we have
  \begin{equation}\label{3.23}
   \int_{\os} |w|^{p+1} \rho |y|^2 \,dy
   \leqslant \frac{2c(p,d_1)}{1-\eta}\ti E_2[w]
         +\frac{\eta}{2(1-\eta)} \int_{\os} w_s^{2} \rho |y|^2 \,dy
         +C(p,d_1,\eta).
  \end{equation}
  Combining \eqref{3.21} with \eqref{3.22} and \eqref{3.23} we obtain that
  \begin{eqnarray*}
   \frac{d(\ti E_2 +c_3 E)}{ds}
   &\leq&
    -2c_4\int_{\os}(|w_s|^2+|\na w|^2)(1+|y|^2) \rho \,dy
    +\frac{2}{p+1}c_3d_2 \int_{\os} |w|^{p+1} \rho \,dy\non\\
   & &
    +\frac{2}{p+1} d_2 \int_{\os} |w|^{p+1} \rho |y|^2 \,dy\non\\
   &\leq&
    \frac{2c(p,d_1)}{1-\eta}\frac{2}{p+1}c_3d_2E[w]
    +\left(\frac{2}{p+1} \frac{\eta}{1-\eta}c_3d_2-2c_4\right)\int_{\os} w_s^{2} \rho \,dy \non\\
   & &
    +C(p,d_1,d_2,\eta)
    +\frac{2c(p,d_1)}{1-\eta}\frac{2}{p+1}c_3d_2\ti E_2[w] \non\\
   & &
    +\left(\frac{2}{p+1} \frac{\eta}{2(1-\eta)}d_2-2c_4\right)\int_{\os} w_s^{2} \rho |y|^2 \,dy
    -c_4\int_{\os} |\na w|^2 (1+|y|^2) \rho \,dy,
  \end{eqnarray*}
  where $d_2$ is a constant such that
  $\sup\left|\dsf{\pa \vv}{\pa s}\right| \leqslant d_2$.
  Take $\eta \leqslant \frac{1}{2}$ small enough such that
  $\dsf{\eta d_2}{(p+1)(1-\eta)} \leqslant \dsf{c_4}{c_3}$, then
  \begin{eqnarray*}
   \frac{d(\ti E_2 +c_3 E)}{ds}
   &\leq&
    -c_4\int_{\os}(|w_s|^2+|\na w|^2)(1+|y|^2) \rho \,dy
    +\frac{8}{(p-1)d_1}c_3d_2E[w]\\
   & &
    +\frac{8}{(p-1)d_1} d_2\ti E_2[w]
    +C(p,d_1,d_2,\eta).
  \end{eqnarray*}
Denote
  $\lambda=\dsf{8}{p-1}\dsf{d_2}{d_1}$, then we get
  $$
  \frac{d(\ti E_2 +c_3 E)}{ds}
  \leq
  -c_4\int_{\os}(w_s^2+|\na w|^2)(1+|y|^2) \rho \,dy
  +\lambda (\ti E_2 +c_3 E) +c_5,
  $$
  where $c_5$ depends on $p, d_1, d_2, \eta$.

\begin{lem}\label{lem11}
 $\ti E_2 +c_3 E \geqslant -\bar{C},$
 where $\bar{C}$ depends on $p, d_1, d_2, \eta$.
\end{lem}
\textbf{Proof.}
 From \eqref{3.14}, using Jensen's inequality, we have
 \begin{eqnarray*}
  \frac{1}{2}\frac{d}{ds}\int_{\os} w^2 \rho |y|^2 \,dy
  & = & -2\ti E_2[w]
        +\frac{p-1}{p+1}\int_{\os} \vv |w|^{p+1} \rho |y|^2 \,dy
        +2\int_{\os} \left(n-\frac{|y|^2}{2}\right)w^2 \rho \,dy \\
  &\geqslant& -2\ti E_2[w]
        -\int_{\os} w^{2} \rho |y|^2 \,dy
        +C\int_{\os} |w|^{p+1} \rho |y|^2 \,dy \\
  &\geqslant& -2\ti E_2[w]
        +(C-\ve)\int_{\os} |w|^{p+1} \rho |y|^2 \,dy
        -C(\ve) \\
  &\geqslant& -2\ti E_2[w]
        -C(\ve)
        +C\left(\int_{\os} w^{2} \rho |y|^2 \,dy\right)^{\frac{p+1}{2}}.
 \end{eqnarray*}
 This inequality plus $c_3 \times$ \eqref{3.9} leads to
 \begin{eqnarray*}
  \frac{1}{2}\frac{d}{ds}\int_{\os} w^2 \rho (|y|^2+c_3) \,dy
  &\geqslant& -2c_3E[w]
        +c_3 C\int_{\os} |w|^{p+1} \rho \,dy
        -2\ti E_2[w]\\
  &  &  +C\left(\int_{\os} w^{2} \rho |y|^2 \,dy\right)^{\frac{p+1}{2}}
        -C(\ve) \\
  &\geqslant& -2\left(\ti E_2 +c_3 E +C(\ve)\right)
        +C\left(\int_{\os} w^{2} \rho (c_3+|y|^2) \,dy\right)^{\frac{p+1}{2}}.
 \end{eqnarray*}
 Denote
 $y(s) \triangleq \int_{\os} w^{2} \rho (c_3+|y|^2) \,dy, \
  J \triangleq \ti E_2 +c_3 E, \
  \bar{C} \triangleq \max\{C(\ve), \dsf{c_5}{\ve}\}$. Then
 \begin{equation}\label{3.25}
  \frac{1}{2}\frac{d}{ds}y(s) \geqslant -2(J+\bar{C}) +Cy^{\frac{p+1}{2}}(s).
 \end{equation}
 We claim that
 $$J \geqslant -\bar{C}.$$
 If not, there exists $s_1$ such that $J(s_1) < -\bar{C}$, then
 \eqref{3.24} tells us that
 $$\left.\frac{d(J+\bar{C})}{ds}\right|_{s_1}
   \leqslant \left.\ve\left(J+\frac{c_5}{\ve}\right)\right|_{s_1}
   \leqslant \ve(J+\bar{C})
   < 0,$$
 which shows that
 $$
 J(s) < -\bar{C} \ \ \forall \ s \geqslant s_1.
 $$
 Therefore from \eqref{3.25} we get
 $\dsf{1}{2}\dsf{d}{ds}y(s) \geqslant Cy^{\frac{p+1}{2}}(s).$
 From this inequality, we easily conclude that $y(s)$ will blow up in
 finite time, which is impossible.
 Hence our lemma holds.

To obtain rough estimates for the higher level energies, the
following two inequalities, i.e. \eqref{1} and \eqref{2}, play an
important role. By Proposition \ref{p3.4} and Young's inequality, we
have
\begin{eqnarray}\label{18}
  \frac{d}{ds}E_{2k}[w]
  &=& -\int_{\os} w_s^2 \rho |y|^{2k}\,dy
      -2k\int_{\os} \rho (y \cdot \na w)w_s |y|^{2k-2}\,dy \non\\
  & &   -\frac{1}{4}\int_{\pa\os} \left|\frac{\pa w}{\pa \gamma}\right|^2
      (y \cdot \gamma)\rho |y|^{2k}\,d\sigma
      -\frac{1}{p+1}\int_{\os} \frac{\pa \vv}{\pa s} |w|^{p+1} \rho |y|^{2k}\,dy\non\\
  &\leqslant& -(1-\ve)\int_{\os} w_s^2 \rho |y|^{2k}\,dy + C(k,\ve)
      \int_{\os} |\na w|^2 \rho |y|^{2k-2}\,dy \non\\
  & & -\frac{1}{p+1}\int_{\os} \frac{\pa \vv}{\pa s} |w|^{p+1} \rho |y|^{2k}\,dy.
\end{eqnarray}
Similar to \eqref{3.23}, we have
\begin{equation}\label{11}
   \int_{\os} |w|^{p+1} \rho |y|^{2k} \,dy
   \leqslant \frac{2c(p,d_1)}{1-\eta} E_{2k}[w]
         +\frac{\eta}{2(1-\eta)} \int_{\os} w_s^{2} \rho |y|^{2k} \,dy
         +C(p,d_1,\eta).
\end{equation}
Taking $\ve, \eta >0$ small enough, we obtain that
\begin{equation}\label{1}
\frac{d}{ds}E_{2k}[w] \leqslant -\frac{1}{2}\int_{\os} w_s^2 \rho
|y|^{2k}\,dy +\mu E_{2k}[w] +C(\mu) + C(\mu)\int_{\os} |\na w|^2
\rho |y|^{2k-2}\,dy,
\end{equation}
for all $\mu \geqslant \lambda$.

On the other hand, by Proposition \ref{p3.3}, H\"{o}lder inequality,
Young's inequality  and Jensen's inequality we have
\begin{eqnarray}\label{2}
\frac{1}{2}\frac{d}{ds}\int_{\os} w^2 |y|^{2k} \rho \,dy
&=& -2E_{2k}[w]+\frac{p-1}{p+1}\int_{\os} \vv |w|^{p+1} \rho |y|^{2k}\,dy \non\\
& & +\int_{\os} k\left(n+2k-2-\frac{1}{2}|y|^2\right)w^2 |y|^{2k-2}\rho\,dy \non\\
&\geqslant & -2E_{2k}[w] - C \int_{\os}w^2 |y|^{2k}\rho\,dy + C
\int_{\os}|w|^{p+1} |y|^{2k}\rho\,dy \non\\
&\geqslant & -2E_{2k}[w] + (C-\ve)
\int_{\os}|w|^{p+1} |y|^{2k}\rho\,dy -C(\ve) \non\\
&\geqslant & -2E_{2k}[w] -C + C \left(\int_{\os}w^2 |y|^{2k}\rho\,dy
\right)^{\frac{p+1}{2}}.
\end{eqnarray}

Now we get following rough estimates
\begin{lem}\label{lem1}
For any $k\in \mathbb{N}$, there exist positive constants $L_k, M_k,
N_k$ and $Q_k$, such that the following estimates hold:
\[
-L_k e^{2\lambda s} \leqslant E_{2k}[w](s) \leqslant M_k e^{2\lambda
s},
\]
\[
\int_0^{\oo}\!\!e^{-2\lambda s}\!\!\int_{\os}|\na w|^2 \rho
|y|^{2k}\, dy \,ds \leqslant N_k,
\]
\[
\int_{\os}w^2\rho |y|^{2k-2}\, dy \leqslant Q_k e^{2\lambda s}
\]
for all $k\in \mathbb{N}$ and $s\geqslant 0$.
\end{lem}
\textbf{Proof.} \ Let $\{\lambda_k\}_{k=1}^{\oo}\subset [\lambda,
2\lambda]$ be a strictly increasing sequence. It suffices to show
the following estimates:
\begin{equation}\label{3}
-L_k e^{\lambda_k s} \leqslant E_{2k}[w](s) \leqslant M_k
e^{\lambda_k s},
\end{equation}
\begin{equation}\label{4}
\int_0^{\oo}\!\!e^{-\lambda_k s}\!\!\int_{\os}|\na w|^2 \rho
|y|^{2k}\, dy \,ds \leqslant N_k,
\end{equation}
\begin{equation}\label{5}
\int_{\os}w^2\rho |y|^{2k-2}\, dy \leqslant Q_k e^{\lambda_k s}.
\end{equation}
We prove these estimates by induction.

\textbf{Step 1.} These estimates holds for $k=1$.

Note that \eqref{3.24} gives us
$\dsf{d}{ds}\left(J+\frac{c_5}{\lambda}\right)
 \leqslant \lambda \left(J+\frac{c_5}{\lambda}\right),$
which imply that $J \leqslant C e^{\lambda s}$. Therefore we now
have $-\bar{C} \leqslant J \leqslant C e^{\lambda s}$ by Lemma
\ref{lem11}.
Using the similar trick of getting \eqref{3.24}, we can write
\eqref{3.13} as a more refinement form:
$$\dsf{d}{ds}\left(E[w]+\frac{c_2}{c_1}\right)
 \leqslant
 \lambda \left(E[w] +\frac{c_2}{c_1}\right),$$
then $E[w] \leqslant C e^{\lambda s}$ and therefore $\ti E_2[w]
\geqslant -\bar{C} -c_3 E[w] \geqslant -C e^{\lambda s}.$ It follows
that
\begin{equation}\label{6}
|\ti E_2[w]|\leqslant C e^{\la s}.
\end{equation}

From \eqref{3.24}, we have
$\dsf{d}{ds}\left(J+\frac{c_5}{\lambda}\right)
 \leqslant
 -c_4\int_{\os}(w_s^2+|\na w|^2)(1+|y|^2) \rho \,dy
  +\lambda \left(J+\frac{c_5}{\lambda}\right).$
Multiplying $e^{-\lambda s}$ on both sides and integrating from $0$
to $\oo$, we obtain that
\begin{equation}\label{3.27}
 \int_0^\oo e^{-\lambda s}\!\!\!\int_{\os}(w_s^2+|\na w|^2)(1+|y|^2) \rho \,dy ds
 \leqslant C.
\end{equation}
In particular, \eqref{4} holds for $k=1$.

Denote $y(s)=\Ios w^2 \rho \,dy$. Notice that
\begin{eqnarray*}
\frac{d}{ds}\int_{\os} w^2 \rho \,dy &=& -4E[w] +
\frac{2(p-1)}{p+1}\Ios \bar V |w|^{p+1}\rho \,dy\\
&\geqslant & -C e^{\la s} +C \left(\Ios w^2 \rho
\,dy\right)^{\frac{p+1}{2}}\\
&=& c_7\left(-c_8e^{\la s} +\left(\Ios w^2 \rho
\,dy\right)^{\frac{p+1}{2}}\right).
\end{eqnarray*}
If there exists $s_1 \geqslant 0$ such that $y(s_1)- 2c_8 e^{\la
s_1}>0$, then at $s_1$,
\begin{eqnarray*}
\left.\frac{d}{ds}(y(s)-2c_8 e^{\la s})\right|_{s_1} &=& y'(s_1)-2\la
c_8 e^{\la s_1}\\
&\geqslant & c_7\left(y(s_1)^{\frac{p+1}{2}}-c_8 e^{\la s_1}\right)
-2\la c_8
e^{\la s_1}\\
&=& c_7\left(y(s_1)^{\frac{p+1}{2}}-c_8(1 +2\la /c_7) e^{\la
s_1}\right)\\
&> & c_7\left(c_8^{\frac{p+1}{2}}e^{\frac{p+1}{2}\la s_1}-c_8(1
+2\la /c_7) e^{\la s_1}\right)\\
&>& 0,
\end{eqnarray*}
since $c_8$ can be large enough. It follows that $y(s) >2 c_8 e^{\la
s}$ for all $s>s_1$. So $y(s)^{\frac{p+1}{2}}> y(s)>2 c_8 e^{\la s}$
and then $\dsf{d}{ds}y(s)\geqslant \frac{c_8}{2}
y^{\frac{p+1}{2}}(s)$ for all $s>s_1$, which implies that $y$ will
blow up in finite time. This contradicts the fact that $y$ is
globally defined. So we have
\begin{equation}\label{7}
y(s)\leqslant 2c_8 e^{\la s}, \ \ \forall \ s\geqslant 0.
\end{equation}
In other words, \eqref{5} holds for $k=1$.

By \eqref{1},
\[
\frac{d}{ds}\left(e^{-\lambda s}E_2[w]\right) \leqslant C
e^{-\lambda s} \int_{\os} |\na w|^2 \rho \,dy + C e^{-\lambda s}.
\]
It follows from \eqref{3.27} that
\[
E_2[w]\leqslant C e^{\lambda s}.
\]

On the other hand, by \eqref{6} and the definition of $\ti E_2$, we
have
\begin{eqnarray*}
-C e^{\la s} &\leqslant & \ti E_2[w] =E_2[w]-\frac{1}{2}\Ios
\left(\frac{1}{2}|y|^2-n\right)w^2\rho \,dy\\
&\leqslant & E_2[w] +\frac{n}{2}\Ios w^2 \rho \,dy\\
&\leqslant & E_2[w] +C e^{\la s},
\end{eqnarray*}
where the last inequality follows from \eqref{5} for $k=1$.
Therefore \eqref{3} also holds for $k=1$.

\textbf{Step 2.} \eqref{3}-\eqref{5} holds for all $k\in
\mathbb{N}$.

Suppose \eqref{3}-\eqref{5} holds for $k\leqslant n$. Since
\eqref{3} holds for $k=n$,  by \eqref{2} and a similar argument to
derive \eqref{7} we conclude that \eqref{5} holds for $k=n+1$. By
\eqref{1}, we have
\[
\frac{d}{ds}(e^{-\la_n s}E_{2n+2}[w])\leqslant  Ce^{-\la_n s}\Ios
|\na w|^2 \rho |y|^{2n} \,dy + C e^{-\la_n s}.
\]
Since \eqref{4} holds for $k=n$, we have
\[
e^{-\la_n s}E_{2n+2}[w] \leqslant C_n.
\]

Now we need to obtain the lower bound for $E_{2n+2}[w]$. Denote
\begin{eqnarray*}
y(s)&=& \Ios w^2 \rho |y|^{2n+2} \,dy\\
z(s)&=& E_{2n+2}[w]+C(\la_n).
\end{eqnarray*}
Then it follows from \eqref{1} and \eqref{2} that
\begin{eqnarray}\label{8}
y'(s)&\geqslant& -4 z(s) + Cy^{\frac{p+1}{2}}(s)\\ \label{9}
z'(s)&\leqslant& \la_n z(s) + C\Ios |\na w|^2 \rho |y|^{2n}\,dy.
\end{eqnarray}
The last inequality implies that
\begin{equation}\label{9}
\frac{d}{ds}(e^{-\la_n s}z(s)) \leqslant e^{-\la_n s} h(s),
\end{equation}
where $h(s)=C\Ios |\na w|^2 \rho |y|^{2n}\,dy$. By induction
hypothesis, we have
\begin{equation}\label{12}
\int_0^{\oo}\!\!e^{-\la_n s}\!\!\Ios |\na w|^2 \rho |y|^{2n}\,dy
\leqslant C_n.
\end{equation}

We claim that
\begin{equation}\label{10}
z(s)\geqslant -N e^{\la_n s}, \ \ \forall \ s\geqslant 0,
\end{equation}
where $N=\int_0^{\oo}e^{-\la_n s}h(s)ds<\oo$.

Otherwise, there exists $s_1\geqslant 0$ such that $e^{-\la_n
s_1}z(s_1)+N<0$. By \eqref{9}, we have
\[
e^{-\la_n s}z(s) - e^{-\la_n s_1}z(s_1)\leqslant \int_{s_1}^s
e^{-\la_n \tau}h(\tau)\,d\tau \leqslant N,
\]
for all $s>s_1$. So $e^{-\la_n s}z(s)\leqslant N + e^{-\la_n
s_1}z(s_1) <0$, i.e., $z(s)<0$ for all $s>s_1$. Now from \eqref{8}
we conclude that $y'(s)\geqslant C y^{\frac{p+1}{2}}(s)$ for all
$s\geqslant s_1$, which implies $y(s)$ blows up in finite time. This
is a contradiction. Therefore $E_{2n+2}[w]\geqslant -C e^{\la_n s}$
and then $|E_{2n+2}[w]|\leqslant C e^{\la_n s}$. In particular,
\eqref{3} holds for $k=n+1$.

Finally, by \eqref{1}, we have
\[
\frac{d}{ds}E_{2n+2}[w]\leqslant -\frac{1}{2}\Ios \! w_s^2\rho
|y|^{2n+2} \,dy + C \Ios \!|\na w|^2\rho |y|^{2n}\,dy +C +\la_n
E_{2n+2}[w].
\]
Combining this with the fact that $|E_{2n+2}[w]|\leqslant C e^{\la_n
s}$ and \eqref{12} we have
\[
\int_0^{\oo}\!\!e^{-\la_n s} \!\!\Ios w_s^2\rho |y|^{2n+2} \,dy \,ds
\leqslant C.
\]
By \eqref{11}, we obtain
\begin{eqnarray*}
\Ios \!|\na w|^2\rho |y|^{2n+2}\,dy &\leqslant & 2E_{2n+2}[w]
+\frac{2}{p+1} \Ios \bar V |w|^{p+1}\rho |y|^{2n+2} \,dy\\
&\leqslant & CE_{2n+2}[w] +C +C \Ios w_s^2\rho |y|^{2n+2} \,dy.
\end{eqnarray*}
Therefore,  by $|E_{2n+2}[w]|\leqslant C e^{\la_n s}$, we get
\begin{eqnarray*}
&&\int_0^{\oo}\! e^{-\la_{n+1} s}|\na w|^2\rho |y|^{2n+2}\,dy\\
&\leqslant& C\int_0^{\oo}\!\! (E_{2n+2}[w]+1)e^{-\la_{n+1} s} ds + C
\int_0^{\oo}\!\!e^{-\la_n s} \!\!\Ios w_s^2\rho |y|^{2n+2}
\,dy \,ds \\
&\leqslant & C\int_0^{\oo}\!\!e^{(\la_n-\la_{n+1}) s} ds +C\\
&\leqslant & C.
\end{eqnarray*}
Hence \eqref{4} holds for $k=n+1$. The Lemma is proved.

\begin{rem}\label{rem}
We have seen in the proof of this Lemma that
\[
-L\leqslant E[w] \leqslant C e^{\la s},
\]
and
\[
\int_0^{\oo}\!\!e^{-\la s}\!\!\int_{\os}|\na w|^2 \rho \, dy \,ds
\leqslant C.
\]
\end{rem}

Next, we need the following
\begin{lem}\label{lem2}
Suppose $\la >\frac{1}{4}$ and for some $\alpha \in (\frac{1}{2},
2\la]$, there exist positive constants $M_k$ and $N_k$, such that
\[
|E_{2k}[w](s)| \leqslant M_k e^{\alpha s},
\]
\[
\int_0^{\oo}\!\!e^{-\alpha s}\!\!\int_{\os}|\na w|^2 \rho |y|^{2k}\,
dy \,ds \leqslant N_k,
\]
hold for all $k\in \mathbb{N}\cup \{0\}$ and $s\geqslant 0$. Then
there exist positive constants $M'_k$ and $N'_k$, such that
\[
|E_{2k}[w](s)| \leqslant M'_k e^{(\alpha -\frac{1}{4}) s},
\]
\[
\int_0^{\oo}\!\!e^{-(\alpha -\frac{1}{4}) s}\!\!\int_{\os}|\na w|^2
\rho |y|^{2k}\, dy \,ds \leqslant N_k,
\]
hold for all $k\in \mathbb{N}\cup \{0\}$ and $s\geqslant 0$. Here we
set $E_0[w]=E[w]$.
\end{lem}
\textbf{Proof.} \ Let $\{\de_k\}_{k=0}^{\oo}\subset [\frac{1}{4},
\frac{1}{3}]$ be a strictly decreasing sequence. It suffices to show
the following estimates:
\begin{equation}\label{13}
|E_{2k}[w](s)| \leqslant M_k e^{(\alpha -\de_k) s},
\end{equation}
\begin{equation}\label{14}
\int_0^{\oo}\!\!e^{-(\alpha -\de_k) s}\!\!\int_{\os}|\na w|^2 \rho
|y|^{2k}\, dy \,ds \leqslant N_k.
\end{equation}
We prove these estimates by induction.

\textbf{Step 1.} These estimates hold for $k=0$.

Recalling \eqref{3.11} we have
\begin{eqnarray}\label{3.26}
 \frac{dE}{ds}
 &\leqslant& -\int_{\os} w_s^2 \rho \,dy
       +\int_{\os} \na V \cdot y e^{-s/2} |w|^{p+1} \rho \,dy \non\\
 &\leqslant& -\int_{\os} w_s^2 \rho \,dy
       +C e^{-s/2} \int_{\os} |y| |w|^{p+1} \rho \,dy \non\\
 &\leqslant& -\int_{\os} w_s^2 \rho \,dy
       +C e^{-s/2} \int_{\os} |y|^2 |w|^{p+1} \rho \,dy
       +C e^{-s/2} \int_{\os} |w|^{p+1} \rho \,dy.
\end{eqnarray}

Also we get
\begin{eqnarray*}
 && e^{- s/2}\int_{\os} |y|^2 |w|^{p+1} \rho \,dy\\
 &\leqslant& C e^{-s/2}\left(\int_{\os} |\na w|^2 |y|^2 \rho \,dy
       +\int_{\os} |w|^{p+1} \rho \,dy
       +C E_2[w]
       +C\right) \\
 &\leqslant& C e^{-s/2}\left(\int_{\os} |\na w|^2 |y|^2 \rho \,dy
       +\int_{\os} |w|^{p+1} \rho \,dy
       +C e^{\alpha s}
       +C\right).
\end{eqnarray*}
By \eqref{3.22} and the assumptions of this Lemma, then we get
\begin{eqnarray}\label{15}
\frac{d}{ds}E[w] &\leqslant & -\frac{1}{2}\Ios w_s^2 \rho \,dy + C
e^{-\frac{s}{2}}\Ios |\na w|^2 \rho |y|^2 \,dy + C e^{(\alpha
-\frac{1}{2})s}+C e^{-\frac{1}{2}s}(E[w]+C)\non\\
&\leqslant &-\frac{1}{2}\Ios w_s^2 \rho \,dy + C
e^{-\frac{s}{2}}\Ios |\na w|^2 \rho |y|^2 \,dy + C e^{(\alpha
-\frac{1}{2})s}.
\end{eqnarray}
So
\[
E[w](s)-E[w](0)\leqslant C\int_0^s\!\!e^{-\frac{\tau}{2}}
\int_{\Om(\tau)}\!\!|\na w|^2\rho |y|^2 \,dy \,d\tau + Ce^{(\alpha
-\frac{1}{2})s}.
\]
We claim that
\begin{equation}\label{16}
\int_0^s e^{-\frac{\tau}{2}}\int_{\Om(\tau)}|\na w|^2 \rho
|y|^{2}\,dy\,d\tau \leqslant C e^{(\alpha -\frac{1}{2})s}.
\end{equation}
Indeed, if we denote the left hand side of \eqref{16} by $f(s)$,
then $\int_0^{\oo}e^{-(\alpha-\frac{1}{2})s}f'(s)\,ds\leqslant C$ by
the assumption. It follows that
\[
C\geqslant \int_0^{s}e^{-(\alpha-\frac{1}{2})s}f'(s)\,ds \geqslant
f(s)e^{-(\alpha-\frac{1}{2})s},
\]
by integration by parts. So \eqref{16} holds and
\[
E[w](s)\leqslant Ce^{(\alpha-\frac{1}{2})s}.
\]
Notice that we have proved that $E[w]\geqslant -L$. Therefore
\eqref{13} holds for $k=0$.

By \eqref{15}, \eqref{16} and $E[w]\geqslant -L$, we deduce that
\begin{equation}\label{17}
\int_0^s \int_{\Om(\tau)}w_s^2 \rho \,dy \,d\tau \leqslant
Ce^{(\alpha-\frac{1}{2})s}.
\end{equation}

As usual, we have
\begin{eqnarray*}
\Ios |\na w|^2 \rho \,dy &\leqslant & 2 E[w] +\frac{2}{p+1}\Ios \bar
V |w|^{p+1}\rho \,dy \\
&\leqslant & C E[w] + C\Ios w_s^2 \rho \,dy +C.
\end{eqnarray*}
Then
\begin{eqnarray*}
e^{-(\alpha -\frac{1}{3}s)}\Ios |\na w|^2 \rho \,dy &\leqslant & C
(E[w]+1)e^{-(\alpha -\frac{1}{3}s)} + C e^{-(\alpha
-\frac{1}{3}s)}\Ios w_s^2 \rho \,dy\\
&\leqslant & C e^{-\frac{1}{6}s}+ C e^{-(\alpha -\frac{1}{3}s)}\Ios
w_s^2 \rho \,dy.
\end{eqnarray*}

Let $f(s)=\displaystyle \int_0^s \int_{\Om(\tau)}w_s^2 \rho \,dy
\,d\tau$. Then for any $s>0$,
\begin{eqnarray*}
\int_0^s \! e^{-(\alpha -\frac{1}{3})\tau}\int_{\Om(\tau)}w_s^2 \rho
\,dy \,d\tau &=& \int_0^s f'(\tau)e^{-(\alpha
-\frac{1}{3})\tau}\,d\tau \\
&=& f(s)e^{-(\alpha -\frac{1}{3})s}+(\alpha -\frac{1}{3})\int_0^s
f(\tau)e^{-(\alpha -\frac{1}{3})\tau}\,d\tau\\
&\leqslant & C,
\end{eqnarray*}
due to \eqref{17}. So
\begin{eqnarray*}
\int_0^{\oo}e^{-(\alpha -\frac{1}{3}\tau)}\int_{\Om(\tau)} |\na w|^2
\rho \,dy \,d\tau &\leqslant & C \int_0^{\oo}
e^{-\frac{1}{6}\tau}\,d\tau + C\int_0^{\oo} e^{-(\alpha
-\frac{1}{3}\tau)}\int_{\Om(\tau)} w_s^2 \rho \,dy\,d\tau\\
&\leqslant & C,
\end{eqnarray*}
i.e., \eqref{14} holds for $k=0$.

\textbf{Step 2.} \eqref{13} and \eqref{14} hold for all $k\in
\mathbb{N}\cup \{0\}$.

Suppose \eqref{13} and \eqref{14} hold for all $k=0, 1, \cdots,
n-1$. Taking $\ve =1/4$ in \eqref{18}, we get
\begin{eqnarray*}
\frac{d E_{2n}[w]}{ds} &\leqslant & -\frac{3}{4}\Ios\!\! w_s^2\rho
|y|^{2n} \,dy +\frac{1}{p+1}\Ios\!\! \left|\frac{\partial \bar
V}{\partial s}\right| |w|^{p+1}\rho |y|^{2n} \,dy + C\Ios\!\! |\na
w|^2
\rho |y|^{2n-2} \,dy\\
&\leqslant & -\frac{3}{4}\Ios\!\! w_s^2\rho |y|^{2n} \,dy
+Ce^{-\frac{s}{2}} \Ios\!\! |w|^{p+1}\rho |y|^{2n+1} \,dy + C\Ios
|\na
w|^2 \rho |y|^{2n-2} \,dy\\
&\leqslant & -\frac{3}{4}\Ios w_s^2\rho |y|^{2n} \,dy + C\Ios |\na
w|^2 \rho |y|^{2n-2} \,dy\\
&& +Ce^{-\frac{s}{2}}\left( \Ios |\na w|^2\rho |y|^{2n+2} \,dy +\Ios
|w|^{p+1}\rho |y|^{2n}\,dy +C -CE_{2n+2}[w] \right)\\
&\leqslant & -\frac{1}{2}\Ios w_s^2\rho |y|^{2n} \,dy + C\Ios |\na
w|^2 \rho |y|^{2n-2} \,dy\\
&&+Ce^{-\frac{s}{2}}\Ios |\na w|^2\rho |y|^{2n+2} \,dy +
Ce^{-\frac{s}{2}}(E_{2n}[w]+C) + C e^{(\alpha -\frac{1}{2})s}\\
&\leqslant & -\frac{1}{2}\Ios w_s^2\rho |y|^{2n} \,dy + C\Ios |\na
w|^2 \rho |y|^{2n-2} \,dy\\
&&+Ce^{-\frac{s}{2}}\Ios |\na w|^2\rho |y|^{2n+2} \,dy + C
e^{(\alpha -\frac{1}{2})s}.
\end{eqnarray*}
Notice that we have used that $\left|\frac{\partial \bar V}{\partial
s}\right|\leqslant C|y|e^{-\frac{s}{2}}$ and the assumptions of the
Lemma. Hence we get
\begin{eqnarray*}
E_{2n}[w](s)-E_{2n}[w](0)&\leqslant& C \int_0^s\!\!
e^{-\frac{\tau}{2}}\!\int_{\Om(\tau)}\!\! |\na w|^2\rho |y|^{2n+2}
\,dy \,d\tau + Ce^{(\alpha -\frac{1}{2})s}\\
&& + C\int_0^s\!\! e^{-\frac{\tau}{2}}\!\int_{\Om(\tau)}\!\! |\na
w|^2\rho |y|^{2n-2} \,dy \,d\tau.
\end{eqnarray*}
Since $\displaystyle\int_0^{\oo}\!\! e^{-\alpha s}\!\Ios\!\! |\na
w|^2\rho |y|^{2n+2} \,dy \,ds \leqslant N_{n+1}$, we get
\[
\int_0^s\!\! e^{-\frac{\tau}{2}}\!\int_{\Om(\tau)}\!\! |\na w|^2\rho
|y|^{2n+2} \,dy \,d\tau\leqslant C e^{(\alpha -\frac{1}{2})s}
\]
as before. Let $f(s)=\displaystyle\int_0^s\!\!\int_{\Om(\tau)}\!\!
|\na w|^2\rho |y|^{2n-2} \,dy \,d\tau$. Then by induction
hypothesis, we have
\[
\int_0^{\oo}\!\! f'(s)e^{-(\alpha -\de_{n-1})s}\,ds\leqslant
N_{n-1}.
\]
So
\begin{eqnarray*}
\int_0^s f'(\tau ) e^{-(\alpha -\de_{n-1})\tau}\,d\tau &=&
f(s)e^{-(\alpha -\de_{n-1})s} +(\alpha -\de_{n-1})\int_0^s f(\tau )
e^{-(\alpha -\de_{n-1})\tau}\,d\tau \\
&\geqslant & f(s)e^{-(\alpha -\de_{n-1})s},
\end{eqnarray*}
i.e., $f(s)\leqslant N_{n-1}e^{(\alpha -\de_{n-1})s}$.

Therefore
\begin{equation}\label{19}
E_{2n}[w]\leqslant N_n e^{(\alpha -\de_{n-1})s}.
\end{equation}

Now let $y(s)=\Ios w^2\rho |y|^{2n}\,dy, \ \ z(s)=E_{2n}[w]+C$. Then
by \eqref{1} and \eqref{2}, we have
\begin{eqnarray*}
y'(s)&\geqslant & -4z(s)+ C y^{\frac{p+1}{2}}(s),\\
z'(s)&\leqslant & 2\la z(s) + C \Ios |\na w|^2 \rho |y|^{2n-2}\,dy
\triangleq  2\la z(s) + h(s).
\end{eqnarray*}
Since $\alpha <2\la, \ z'(s)\leqslant (\alpha -\de'_n)z(s)+g(s)$,
where $g(s)=(2\la -\alpha +\de'_n)z(s)+h(s)$ and $\de'_n \in (\de_n,
\de_{n-1})$. It follows from \eqref{19} and induction hypothesis
that
\begin{eqnarray*}
\int_0^{\oo}e^{-(\alpha -\de'_n)s}g(s)\,ds &\leqslant &
C\!\int_0^{\oo}\!e^{(\de'_n -\de_{n-1})s}\,ds +C\!
\int_0^{\oo}\!e^{-(\alpha -\de'_n)s}\!\Ios \!|\na w|^2 \rho
|y|^{2n-2}\,dy\,ds\\
&\leqslant & C.
\end{eqnarray*}
A similar argument to obtain \eqref{10} gives us
\begin{equation}\label{20}
z(s)\geqslant -Ce^{(\alpha -\de'_n)s}.
\end{equation}
From \eqref{19} and \eqref{20}, we know that \eqref{13} holds for
$k=n$.

From the fact that
\[
\frac{d E_{2n}[w]}{ds} \leqslant -\frac{1}{2}\Ios w_s^2\rho |y|^{2n}
\,dy + (\alpha -\de'_n)E_{2n}[w]+g(s)+C
\]
and above estimates, we have
\[
\int_0^{\oo}e^{-(\alpha -\de'_n)s}\Ios w_s^2\rho |y|^{2n} \,dy\,ds
\leqslant C.
\]

As before, we have
\[
\Ios \!|\na w|^2 \rho |y|^{2n}\,dy\leqslant CE_{2n}[w] +C\Ios \!
w^2_s \rho |y|^{2n}\,dy +C.
\]
Multiplying $e^{-(\alpha -\de_n)s}$ on both sides and integrating
over $(0,\oo)$, we obtain
\begin{eqnarray*}
&&\int_0^{\oo}\!\!e^{-(\alpha -\de_n)s}\!\!\Ios \!|\na w|^2 \rho
|y|^{2n}\,dy \,ds \\
&\leqslant & C \int_0^{\oo}\!\!e^{-(\alpha -\de_n)s}e^{(\alpha
-\de'_n)s}\,ds +C + C\int_0^{\oo}\!\!e^{-(\alpha -\de'_n)s}\!\!\Ios
\! w_s^2
\rho |y|^{2n}\,dy \,ds\\
&\leqslant & C,
\end{eqnarray*}
i.e., \eqref{14} holds for $k=n$. So the proof of this Lemma is
complete.

To obtain the upper bound of $E[w]$, we also need the following
\begin{lem}\label{lem3}
Suppose that there exist two positive constants $M,N$ and some
$\alpha \in (0, \frac{1}{2})$ such that
\[
|E_{2}[w](s)| \leqslant M e^{\alpha s},
\]
\[
\int_0^{\oo}\!\!e^{-\alpha s}\!\!\int_{\os}|\na w|^2 \rho |y|^{2}\,
dy \,ds \leqslant N.
\]
Then we have $$E[w] \leqslant K_2.$$
\end{lem}

\textbf{Proof.} \ Recall from \eqref{3.26} that
\[
 \frac{dE}{ds}
 \leqslant -\int_{\os} w_s^2 \rho \,dy
       +C e^{-s/2} \int_{\os} |y|^2 |w|^{p+1} \rho \,dy
       +C e^{-s/2} \int_{\os} |w|^{p+1} \rho \,dy.
\]

By the lower bound of $E_2$ and Young's inequality, we get
\begin{eqnarray}\label{3.28}
 e^{- s/2}\int_{\os} |y|^2 |w|^{p+1} \rho \,dy
 &\leqslant& C e^{-s/2}\left(\int_{\os} |\na w|^2 |y|^2 \rho \,dy
       +\int_{\os} |w|^{p+1} \rho \,dy
       +C e^{\alpha s}
       +C\right) \non\\
 &\leqslant& C e^{-s/2} \int_{\os} |\na w|^2 |y|^2 \rho \,dy
       +C e^{-s/2} \int_{\os} |w|^{p+1} \rho \,dy \non\\
 &   & +C e^{-s/2}
       +C e^{(\alpha-\frac{1}{2})s}.
\end{eqnarray}
Using \eqref{3.12}, we have
\begin{eqnarray}\label{3.29}
 \frac{dE}{ds}
 &\leqslant& -\int_{\os} w_s^2 \rho \,dy
       +C e^{-s/2} \int_{\os} |\na w|^2 |y|^2 \rho \,dy \non\\
 &   & +C e^{-s/2} \int_{\os} |w|^{p+1} \rho \,dy
       +C e^{-s/2} +C e^{(\alpha-\frac{1}{2})s}\non\\
 &\leqslant& -\frac{1}{2}\int_{\os} w_s^2 \rho \,dy
       +C e^{-s/2} \int_{\os} |\na w|^2 |y|^2 \rho \,dy \non\\
 &   & +C e^{-s/2} (E[w]+C)
       +C e^{(\alpha-\frac{1}{2})s}.
\end{eqnarray}
By Lemma \ref{lem}, we may assume $E[w]+C>1$. So
\[
\frac{d}{ds}\log (E[w]+C)\leqslant C e^{-s/2} \int_{\os} |\na w|^2
|y|^2 \rho \,dy +C e^{-s/2} +C e^{(\alpha-\frac{1}{2})s}.
\]
Noticing that $\alpha <\frac{1}{2}$, we obtain that $E[w] \leqslant
K_2$ from the assumptions.

\textbf{Proof of Proposition \ref{p3.1}} Combining Lemma \ref{lem3}
with Lemma \ref{lem1}, Lemma \ref{lem2} and Remark \ref{rem}, we get
the upper bound of $E[w]$ immediately. Notice that the lower bound
of $E[w]$ has been obtained in Lemma \ref{lem}. So the proof is
complete.

\subsubsection{Proof of Proposition \ref{p3.2}}

\textbf{Proof of \eqref{3.6}.} From \eqref{3.12} we have
$$
\int_{\os} |w|^{p+1} \rho \,dy \leqslant \ve \int_{\os} w_s^2 \rho
\,dy
    +C(\ve).
$$
Then \eqref{3.29} tells us that
\begin{eqnarray*}
 \frac{dE}{ds}
 \leqslant
 \left(-\frac{1}{2}-\ve e^{-s/2}\right) \int_{\os} w_s^2 \rho \,dy
 +C(\ve) e^{-s/2} +f(s),
\end{eqnarray*}
where $f(s)=C e^{-s/2} \int_{\os} |\na w|^2 |y|^2 \rho \,dy$, which
is an integrable function. Integrating this inequality from $s_0$ to
$T$, we get
$$
\frac{1}{4} \int_{s_0}^T \!\!\!\int_{\os} w_s^2 \rho \,dy \leqslant
\int_{s_0}^T \left(C e^{-s/2} +f(s)\right)\,ds
    +E(s_0)
    -E(T).
$$
It follows that
$$\int_0^\oo \|w_s;L_\rho^2(\os)\|^2 ds \leqslant N_1.$$

\textbf{Proof of \eqref{3.7}.} Making use of Jensen's inequality,
from \eqref{3.9}, we get
$$
\frac{1}{2} \frac{d}{ds} \int_{\os} w^2 \rho \,dy \geqslant
-2K_2+C(p,d_2,\Om)\left(\int_{\os} w^2 \rho
\,dy\right)^{\frac{p+1}{2}}.
$$
We assert that
$$
\int_{\os} w^2 \rho \,dy \leqslant N_2,
$$
where $N_2=\left(\dsf{2K_2}{C(p,d_2,\Om)}\right)^{\frac{2}{p+1}}$ is
the zero of $-2K_2+C(p,d_2,\Om)x^{\frac{p+1}{2}}=0$.

If not, there exists $s_1$ such that
$$\int_{\Om(s_1)} w^2 \rho \,dy > \left(\frac{2K_2}{C(p,d_2,\Om)}\right)^{\frac{2}{p+1}}.$$
Then
$$\left.\frac{1}{2} \frac{d}{ds} \int_{\os} w^2 \rho \,dy\right|_{s=s_1} > C>0,$$
which implies that
$$
\int_{\os} w^2 \rho \,dy > 2C \ \ \forall \ s>s_1.
$$
Then there exists some $\bar{t}$ such that for $s > \bar{t}$,
$$-2K_2+C(p,d_2,\Om)\left(\int_{\os} w^2 \rho \,dy\right)^{\frac{p+1}{2}}
\geqslant \frac{C(p,d_2,\Om)}{2}\left(\int_{\os} w^2 \rho
\,dy\right)^{\frac{p+1}{2}}$$ so that $y$ blows up in finite time,
which is impossible.

\textbf{Proof of \eqref{3.8}.} Recall that $\vv \geqslant d_1$ and
$E[w] \leqslant K_2$. Then from \eqref{3.9} we see that
$$
\int_{\os} |w|^{p+1} \rho \,dy \leqslant \ve
\frac{2(p+1)}{d_1(p-1)}K_2
    +\frac{p+1}{d_1(p-1)}\left(\int_{\os} |w|^2 \rho \,dy\right)^{\frac{1}{2}}
     \left(\int_{\os} |w_s|^2 \rho \,dy\right)^{\frac{1}{2}}.
$$
Therefore by \eqref{3.6} and \eqref{3.7} we have
$$
\int_{s}^{s+1}\!\!\!\left(\int_{\os} |w|^{p+1} \rho \,dy\right)^2ds
\leqslant C +CN_2\int_{0}^{\oo}\!\!\!\int_{\os} |w_s|^2 \rho \,dy
\leqslant N_3.
$$

\subsection{Proof of Theorem \ref{thm1}}

Let $\psi \in C^2(\mathbb{R}^n)$ be a bounded function with
$\mbox{supp} \psi \subset B_{2R}(0) \cap \Om.$ Then $\psi w$
satisfies
\begin{eqnarray}\label{5.1}
 \rho (\psi w)_s
 -\na \cdot \left(\rho \na (\psi w)\right)
 +\di (\rho w \na \psi) +\rho \na \psi \cdot \na w
  + \beta \psi \rho w -\vv \psi |w|^{p-1}w \rho=0 \non\\
  \ \ \ \ \ \ \mbox{\rm in} \ \ \os \times (0, \oo).
\end{eqnarray}
We introduce two types of local energy.
\begin{eqnarray}
 && \hspace{-1cm}E_{\psi}[w](s)=\dsf{1}{2}\int_{\os}\!\! \left(|\na (\psi w)|^2
     +(\beta \psi^2-\na |\psi|^2) w^2\right)\rho \,dy
     -\dsf{1}{p+1}\int_{\os} \!\!\vv \psi^2 |w|^{p+1}\rho \,dy, \label{5.2}\\
 && \hspace{-1cm}\mathcal{E}_{\psi}[w](s)=\dsf{1}{2}\int_{\os} \psi^2(|\na w|^2 +\beta w^2)\rho \,dy
     -\dsf{1}{p+1}\int_{\os} \vv \psi^2 |w|^{p+1}\rho \,dy. \label{5.3}
\end{eqnarray}
By the similar trick of \cite{gms1}, we could establish a lower and
an upper bound for $\mathcal{E}_{\psi}[w]$. We just list some
important results and ignore the proof.

\subsubsection{Upper bound for $\mathcal{E}_{\psi}[w]$}

Using \eqref{3.5} and \eqref{3.7} we obtain that
\begin{equation}\label{5.4}
 \|w(s);W_\rho^{1,2}(\os)\|^2
 \leqslant K_1(1 +\|w_s(s);L_\rho^2(\os)\|)
 \ \ \ \mbox{\rm for all} \ s \geqslant 0,
\end{equation}
where $\|w(s);W_\rho^{1,2}\left(\os\right)\|^2
=\beta\|w(s);L_\rho^2(\os)\|^2 +\|\na w(s);L_\rho^2(\os)\|^2.$

\begin{prop}(Quasi-monotonicity of $\mathcal{E}_{\psi}[w]$)
\begin{equation}\label{5.5}
 \frac{d}{ds}\mathcal{E}_{\psi}[w](s)
 \leqslant
 L_1(1 +\|w_s(s);L_\rho^2(\os)\|)
 +C e^{-s/2} \int_{\os} \psi^2 |y| |w|^{p+1} \rho \,dy
\end{equation}
\end{prop}
for all $s>0.$

\begin{prop} There exists a positive constant $K_2$, such that
\begin{equation}\label{5.6}
 \int_s^{s+1}\mathcal{E}_{\psi}[w](\tau)\,d \tau \leqslant K_2
 \ \ \ \mbox{for all} \ s \geqslant 0,
\end{equation}
where $K_2$ depends on $n, p, \|\psi\|_{\oo}$, upper bound for
$\mathcal{E}_{\psi}[w]$ and upper bound for $\vv$.
\end{prop}

Note that
$$\int_s^{s+1}C e^{-\tau/2} \!\!\int_{\Om(\tau)} \psi^2 |y| |w|^{p+1} \rho \,dy d\tau
 \leqslant C.$$
Thanks to \eqref{5.5}, \eqref{3.6} and \eqref{5.6} we can derive an
upper bound for $\mathcal{E}_{\psi}[w].$

\begin{thm}
\begin{equation}\label{5.7}
 \mathcal{E}_{\psi}[w] \leqslant M
 \ \ \ \mbox{for all} \ s \geqslant 0.
\end{equation}
\end{thm}

\subsubsection{Lower bound for $\mathcal{E}_{\psi}[w]$}

Notice that
$$E_\psi-\mathcal{E}_{\psi}
 =\int_{\os} \psi w (\na \psi \cdot \na w) \rho \,dy.$$

By estimating $|E_\psi-\mathcal{E}_{\psi}|$ and using \eqref{3.7} we
obtain

\begin{prop}
 There exists a positive constant $J_1$ such that
 \begin{equation}\label{5.8}
  \frac{1}{2} \frac{d}{ds} \int_{\os} |\psi w|^2 \rho \,dy
  \geqslant
  -2\mathcal{E}_{\psi}
  -J_1
  +\frac{p-1}{p+1}\int_{\os} \vv \psi^2 |w|^{p+1} \rho |y|^2 \,dy.
 \end{equation}
\end{prop}

By \eqref{5.8}, \eqref{5.5} and \eqref{3.6} we obtain that
\begin{thm}
 There exists a positive constant $L_2$ such that
 \begin{equation}\label{5.9}
  \mathcal{E}_{\psi}[w](s)
  \geqslant
  -L_2
  \ \ \ \mbox{for all} \ s\geqslant0.
 \end{equation}
\end{thm}

Once we have these bounds for the local energies, the proof of
Theorem \ref{thm1} follows from bootstrap arguments, an
interpolation theorem in \cite{cl} and the interior regular theorem
in \cite{lsu} as in \cite{gms1,gms2}. We omit the details since
there is no anything new.

\begin{rem}
If we only treat nonnegative solution to \eqref{1.1}, then Theorem
\ref{thm1} can be proved through the bounds we have obtained in
Section 2.1. We can combine the methods in \cite{gk2} and \cite{pqs}
to get the blow-up rate estimate.
\end{rem}

\section{Asymptotic behavior of the Blow-Up Time and Blow-Up set}

In this section, we are interested in the following problem
\begin{equation*}
\left\{
\begin{array}{ll}
u_t=\Delta u+V(x) u^{p}  & \ \ \mbox{in} \ \ \Omega \times
(0,T),\\
u(x,t)= 0 & \ \ \mbox{on}  \ \partial \Omega \times
(0,T),\\
u(x,0)=M\vf(x) & \ \ \mbox{in} \ \ \Omega,
\end{array}
\right.
\end{equation*}
where $\vf \in C(\bar \Om)$ satisfies $\vf|_{\partial \Om} =0,
\vf(x)>0, \ \forall \ x\in \Om$ and $V$ satisfies the conditions
described as in Section 1.

The main goal of this section is to prove Theorem \ref{thm2} and
\ref{thm3}.

\textbf{Proof of Theorem \ref{thm2}.} \ That blow-up occurs for
large $M$ is standard fact.  Let $\bar{a} \in \Om$ such that
$\vf^{p-1}(\bar{a})V(\bar{a})=\max\limits_{x}\vf^{p-1}(x)V(x).$

Since $\vf$ and $V$ are continuous, it follows that $\forall \
\ve>0, \exists \ \de>0,$ such that
$$V(x)>V(\bar{a})-\dsf{\ve}{2}, \ \ \vf(x)>\vf(\bar{a})-\dsf{\ve}{2},
 \ \ \forall \ x \in B(\bar{a},\de).$$

Let $w$ be the solution of
\begin{equation}\label{2.1}
\left\{
 \begin{array}{ll}
  w_t=\Delta w+\left(V(\bar{a})-\dsf{\ve}{2}\right)w^p
   & \ \ \mbox{in} \ \ B(\bar{a},\de) \times (0,T_w),\\
  w= 0 & \ \ \mbox{on}  \ \pa B(\bar{a},\de) \times (0,T_w),\\
  w(x,0)=M\left(\vf(\bar{a})-\ve\right) & \ \ \mbox{in} \ \ B(\bar{a},\de)
 \end{array}
\right.
\end{equation}
and $T_w$ its corresponding blow up time.

A comparison argument shows that $u \geq w$ in $B(\bar{a},\de)
\times (0,T)$ and hence $T \leq T_w.$

Our goal is to estimate $T_w$ for large values of $M$. Define
$$I(w)=\dsf{1}{2}\int_{B(\bar{a},\de)} |\na w|^2 \,dx
     -\frac{V(\bar{a})-\dsf{\ve}{2}}{p+1}\int_{B(\bar{a},\de)} w^{p+1} \,dx,$$
then
\begin{eqnarray*}
  I'(t)
  &=&\int_{B(\bar{a},\de)} \na w \cdot \na w_t \,dx
     -\left(V(\bar{a})-\dsf{\ve}{2}\right)\int_{B(\bar{a},\de)} w^pw_t \,dx\\
  &=&-\int_{B(\bar{a},\de)}w_t\left(\De w+\left(V(\bar{a})-\dsf{\ve}{2}\right)w^p\right) \,dx\\
  &=&-\int_{B(\bar{a},\de)} w_t^2 \,dx.
\end{eqnarray*}
Set $\Phi(t)=\dsf{1}{2}\int_{B(\bar{a},\de)} w^2(x,t)\,dx,$ then we
obtain that
\begin{eqnarray}\label{2.2}
 \Phi'(t) &=& \int_{B(\bar{a},\de)} ww_t \,dx\non\\
 &=& \int_{B(\bar{a},\de)} w\left(\Delta w+\left(V(\bar{a})-\dsf{\ve}{2}\right)w^p\right) \,dx\non\\
 &=& -\int_{B(\bar{a},\de)} |\na w|^2
     +\left(V(\bar{a})-\dsf{\ve}{2}\right)\int_{B(\bar{a},\de)} w^{p+1} \,dx\non\\
 &=& -2I(w)
     +\frac{p-1}{p+1}\left(V(\bar{a})-\dsf{\ve}{2}\right)\int_{B(\bar{a},\de)} w^{p+1} \,dx\non\\
 &>& -2I(w)
       +\dsf{p-1}{p+1}\left(V(\bar{a})- \ve \right)|B|^{\frac{1-p}{2}}
        \left(\int_{B(\bar{a},\de)} w^{2} \,dx\right)^{\frac{1+p}{2}}\non\\
 &=& -2I(w_0)
     +2\int_0^t\!\!\!\int_{B(\bar{a},\de)} w_t^2 \,dxdt
     +\ti{C}\Phi^{\frac{1+p}{2}}(t),
\end{eqnarray}
where $\ti{C}=\dsf{p-1}{p+1}\left(V(\bar{a})- \ve
\right)|B|^{\frac{1-p}{2}}2^{\frac{1-p}{2}}.$

In particular, $\Phi'(t)>0.$

On the other hand,
$$\Phi'(t) =  \int_{B(\bar{a},\de)} ww_t \,dx
          \leqslant \left(\int_{B(\bar{a},\de)} w^{2} \,dx\right)^{\frac{1}{2}}
               \left(\int_{B(\bar{a},\de)} w_t^{2} \,dx\right)^{\frac{1}{2}}
           =  \left(2\Phi(t)\right)^{\frac{1}{2}}
               \left(\int_{B(\bar{a},\de)} w_t^{2} \,dx\right)^{\frac{1}{2}},$$
which tells us that $\displaystyle\int_{B(\bar{a},\de)} w_t^{2} \,dx
 \geqslant \dsf{(\Phi'(t))^2}{2\Phi(t)}.$
Therefore from \eqref{2.2} we get
$$\Phi'(t)
  >-2I(w_0)
   +\int_0^t\dsf{(\Phi'(t))^2}{\Phi(t)}\,dt
   +\ti{C}\Phi^{\frac{1+p}{2}}(t).$$

Set
$f(t)=-2I(w_0)+\displaystyle\int_0^t\dsf{(\Phi'(t))^2}{\Phi(t)}\,dt$
and
     $g(t)=\dsf{2}{p-1}\ti{C}\Phi^{\frac{1+p}{2}}(t).$

Note that
\begin{eqnarray*}
 &&f(0)=-2I(w_0)
     =\frac{2}{p+1}\left(V(\bar{a})-\dsf{\ve}{2}\right)|B|M^{p+1}(\vf(\bar{a})-\ve)^{p+1},\\
 &&g(0)=\frac{2}{p+1}(V(\bar{a})-\ve)|B|M^{p+1}(\vf(\bar{a})-\ve)^{p+1}.
\end{eqnarray*}
It follows that $f(0)>g(0).$ Hence
$$\Phi'(0)
 > f(0)+\ti{C}\Phi^{\frac{1+p}{2}}(0)
 > g(0)+\ti{C}\Phi^{\frac{1+p}{2}}(0)
 = \frac{p+1}{p-1}\ti{C}\Phi^{\frac{1+p}{2}}(0).$$
Then $\exists \ \eta>0,$ such that $\Phi'(t)
 \geqslant \dsf{p+1}{p-1}\ti{C}\Phi^{\frac{1+p}{2}}(t),
 \ \ t\in [0,\eta].$

Define $A=\{\theta \in [0, T_\Phi]:
  \Phi'(t) \geqslant \dsf{p+1}{p-1}\ti{C}\Phi^{\frac{1+p}{2}}(t),
  \ \ t\in [0,\theta]\},$
where $T_\Phi$ is the blow-up time of $\Phi.$ Then $A$ is closed. On
the other hand, $A$ is open. In fact, $\forall \ \theta \in A,$
since
$$f'(t)=\dsf{(\Phi'(t))^2}{\Phi(t)}, \ \
  g'(t)=\dsf{p+1}{p-1}\ti{C}\Phi^{\frac{p-1}{2}}(t)\Phi'(t),$$
it follows that $f'(t)>g'(t)$ for $t \in [0, \theta].$

Recall that $f(0)>g(0).$ We conclude that
$$f(t)>g(t), \ \  t\in [0, \theta].$$
In particular, $f(\theta) > g(\theta).$

Thus, there exists $\bar{\beta}>0$ such that for all $\beta \in [0,
\bar{\beta}], \ f(\theta+\beta)>g(\theta+\beta)$ or
$$\Phi'(\theta+\beta)>\dsf{p+1}{p-1}\ti{C}\Phi^{\frac{1+p}{2}}(\theta+\beta),$$
which means $\theta+\bar{\beta} \in A.$ Therefore $A =[0,
T_{\Phi}].$ In other words,
$$\Phi'(t)\geqslant\dsf{p+1}{p-1}\ti{C}\Phi^{\frac{1+p}{2}}(t), \ \ t \in [0, T_{\Phi}].$$
Integrating this inequality from 0 to $T_{\Phi}$, we get
$$T_\Phi \leqslant \frac{1}{(p-1)(V(\bar{a})-\ve)M^{p-1}(\vf(\bar{a})-\ve)^{p-1}}.$$

Since $\ve >0$ is arbitrarily small, the Theorem follows readily
from the above estimate.

\textbf{Proof of Theorem \ref{thm3}.} \ The proof is almost the same
as in \cite{cer}. The only different thing is that we improve their
Lemma 2.2. For the reader's convenience, we outline the proof here.

Let $M$ be large such that the solution $u$ blows up in finite time
$T=T(M)$ and let $a=a(M)$ be a blow-up point. To involve the
information of $T$, we modify the definition of $w$ to be
\[
w(y,s)=(T-t)^{\frac{1}{p-1}}u(a+y(T-t)^{\frac{1}{2}},t)|_{t=T(1-e^{-s})}.
\]
Then $w$ satisfies
\[
\rho w_s=\na \cdot (\rho \na w) -\beta \rho w
+V(a+yT^{\frac{1}{2}}e^{-\frac{s}{2}}) |w|^{p-1}w \rho
  \ \ \ \mbox{\rm in} \ \os \times (0, \oo),
\]
where $\os=\{y | a+yT^{\frac{1}{2}}e^{-\frac{s}{2}} \in \Om\}$.

Consider the frozen energy
\[
E(w)=\Ios \left(\frac{1}{2}|\na w|^2 +\frac{\beta}{2}w^2
-\frac{1}{p+1}V(a)w^{p+1} \right)\rho \,dy.
\]
Then
\begin{eqnarray*}
\frac{dE}{ds}&\leqslant &-\Ios w^2_s \rho \,dy +\Ios
(V(a+yT^{\frac{1}{2}}e^{-\frac{s}{2}})-V(a))w^pw_s \rho \,dy\\
&\leqslant & -\Ios w^2_s \rho \,dy +
CT^{\frac{1}{2}}e^{-\frac{s}{2}} \left(\Ios w^2_s \rho
\,dy\right)^{\frac{1}{2}}.
\end{eqnarray*}
We have used Theorem \ref{thm1} and H\"{o}lder inequality in the
last inequality. So $ \dsf{dE}{ds} \leqslant C T e^{-s},$ and then
\[
E(w)\leqslant E(w_0) +C T.
\]

Since $w$ is bounded, by the argument of \cite{gk2} and \cite{gk3},
we conclude that
\[
\lim_{s\to \oo}w(y,s)=k(a)\triangleq
\frac{1}{((p-1)V(a))^{\frac{1}{p-1}}}
\]
uniformly in any compact set, and
\[
E(w(\cdot, s)) \to E(k(a)) \ \ \ \mbox{as } s\to \oo.
\]
So
\begin{equation}\label{21}
E(k(a))\leqslant E(w_0)+CT.
\end{equation}
By Theorem \ref{thm2}, we estimate $E(w_0)$ to get $ E(w_0)\leqslant
E(T^{\frac{1}{p-1}}M\vf (a)) + CT^{\frac{1}{2}}$. So
\[
E(k(a))\leqslant E(T^{\frac{1}{p-1}}M\vf (a)) + CT^{\frac{1}{2}}
\]
Observe that $E(b)=\Gamma F(b)$ for any constant $b$, where $\Gamma
=\int \rho \,dy$ and
$F(x)=\frac{1}{2\beta}x^2-\frac{1}{p+1}V(a)x^{p+1}$. It follows that
$F$ attains a unique maximum at $k(a)$ and there exist $\alpha,
\beta$ such that if $|x-k(a)|<\alpha$ then $F''(x)<-1/2$ and if
$|F(x)-F(k(a))|<\beta$ then $|x-k(a)|<\alpha$. From \eqref{21}, we
have $F(k(a))\leqslant F(T^{\frac{1}{p-1}}M\vf (a))
+CT^{\frac{1}{2}}$. By the properties of $F$ we have
\[
CT^{\frac{1}{2}} \geqslant F(k(a))- F(T^{\frac{1}{p-1}}M\vf (a))
\geqslant \frac{1}{4}(k(a)- T^{\frac{1}{p-1}}M\vf (a))^2.
\]
By Theorem \ref{thm2}, for any $k>0$ there exists $M_k>0$ such that
if $M>M_k$, we have
\begin{eqnarray*}
k(a)-CT^{\frac{1}{4}}&\leqslant & T^{\frac{1}{p-1}}M\vf (a)\\
&\leqslant & k(a)\theta (a) +\frac{C\vf (a)}{M^k},
\end{eqnarray*}
where
\[
\theta (a)=\frac{\vf (a) V(a)^{\frac{1}{p-1}}}{\vf (\bar a) V(\bar
a)^{\frac{1}{p-1}}},\ \ \ \vf (\bar a) V(\bar a)^{\frac{1}{p-1}}=
\max_{x\in \Om}\vf (x) V(x)^{\frac{1}{p-1}}.
\]
Therefore, we get
\[
k(a)(1-\theta (a))\leqslant \frac{C\vf
(a)}{M^k}+\frac{C}{M^{\frac{p-1}{4}}}\leqslant
\frac{C}{M^{\frac{p-1}{4}}}
\]
if we choose $k > \frac{p-1}{4}$. Then
\[
\theta (a)\geqslant 1-\frac{C}{M^{\frac{p-1}{4}}}.
\]
This implies
\[
\vf (a) V(a)^{\frac{1}{p-1}} \geqslant \vf (\bar a) V(\bar
a)^{\frac{1}{p-1}} -\frac{C}{M^{\frac{p-1}{4}}}.
\]

We can deduce from this inequality that $\vf(a)\geqslant C>0$ for
large $M$. So
\[
\frac{1}{\vf(a)((p-1)V(a))^{\frac{1}{p-1}}}-
\frac{CT^{\frac{1}{4}}}{\vf(a)} \leqslant M T^{\frac{1}{p-1}}.
\]
Therefore
\[
\frac{1}{\vf(\bar a)((p-1)V(\bar a))^{\frac{1}{p-1}}}-
CT^{\frac{1}{4}} \leqslant M T^{\frac{1}{p-1}},
\]
i.e.,
\[
\frac{1}{\vf(\bar a)((p-1)V(\bar a))^{\frac{1}{p-1}}}-
\frac{C}{M^{\frac{1}{p-1}}} \leqslant M T^{\frac{1}{p-1}}.
\]
The Theorem is proved.

\end{document}